\documentclass[11pt,a4paper]{article}

\usepackage[french]{babel}

\usepackage[latin1]{inputenc}

\usepackage[centertags]{amsmath}

\usepackage{amsfonts}

\usepackage{amssymb}

\usepackage{amsthm}

\usepackage{newlfont}

\newlength{\defbaselineskip}

\usepackage{fancyhdr}
 \fancypagestyle{plain}{\fancyhead{}

}

\newtheorem{theo}{Théorème}

\newtheorem{prop}{Proposition}

\newtheorem{defi}{Définition}

\newtheorem{lem}{Lemme}

\newtheorem{coro}{Corollaire}

\newcommand{\complex}{\mathbb C}

\newcommand{\integer}{\mathbb N}

\newcommand{\rinteger}{\mathbb Z}

\newcommand{\real}{\mathbb R}

\newcommand{\puiss}[3]{\underset{#1}{}[#3]\underset{#2}{}}

\newcommand{\deform}[3]{ #3 ^{(#1,#2)}}

\newcommand{\opmixte}[2]{\vphantom{}_{#1}\delta_{#2}}

\newcommand{\la}{\lambda}

\newcommand{\al}{\alpha}

\newcommand{\ve}{\varepsilon}

\newcommand{\resp}{{\sl resp.} }

\newcommand{\dis}{\displaystyle }

\newcommand{\qhfact}[1]{\mathcal{O}_f^{(#1)}}

\begin{document}

\title{\Huge Familles fuchsiennes d'\'equations aux  ($q$-)différences et confluence.\\}

\author{Anne DUVAL et Julien ROQUES.}

\maketitle  

\hrule
 \abstract{On commence par pr\'esenter une m\'ethode de r\'esolution d'une famille de syst\`emes fuchsiens  d'op\'erateurs de pseudo-d\'erivations associ\'ees \`a une famille \`a deux param\`etres d'homographies qui unifie et g\'en\'eralise les cas connus des systèmes différentiels, aux différences ou aux $q$-différences. Nous traitons ensuite dans cette famille des probl\`emes de confluence que l'on peut voir comme des probl\`emes de continuit\'e en ces deux param\`etres.}\\

\hrule
\tableofcontents
$_{}$\linebreak

\hrule

\section{Introduction}

\subsection{Probl\`ematique}

L'\'etude locale des trois types de syst\`emes lin\'eaires suivants :
\begin{enumerate}
\item diff\'erentiel : $Y'(x)=A(x)Y(x)$,
\item aux $q$-diff\'erences : $Y(qx)=A(x)Y(x)$, o\`u $q$ est un nombre complexe non nul de module diff\'erent de $1$,
\item aux diff\'erences : $Y(x-1)=A(x) Y(x)$,
\end{enumerate}
o\`u $x$ est une variable complexe et $A(x)$ une matrice de fonctions, fait appara\^\i tre de grandes similitudes.  Nous nous restreindrons dans ce qui suit à l'\'etude locale de ces systèmes au voisinage de l'infini qui est le seul point fixe, sur la sph\`ere de Riemann, de la translation $x\mapsto x-1$ et l'un des deux points fixes de l'homographie $x\mapsto qx$. 

Une matrice $T(x)$ inversible au voisinage de l'infini d\'efinit une \og transformation de jauge\fg{}, c'est-\`a-dire associe \`a la matrice $A(x)$ la matrice $B(x)$ du syst\`eme que v\'erifie $X(x)=T(x)Y(x)$ lorsque $Y(x)$ v\'erifie le syst\`eme de matrice $A(x)$. On a : \begin{enumerate}
\item $B(x)=T(x)^{-1}A(x)T(x)-T(x)^{-1}T'(x)$ dans le cas diff\'erentiel,
\item $B(x)=T(qx)^{-1}A(x)T(x)$ dans le cas des $q$-diff\'erences,
\item $B(x)=T(x-1)^{-1}A(x)T(x)$ dans le cas des diff\'erences.
\end{enumerate}
Enon\c cons quelques r\'esultats sous une forme faisant ressortir le parall\'elisme.
\begin{enumerate}
\item Supposons que $\displaystyle A(x)$ est la somme d'une s\'erie $\sum_{s\geq 0}A_s x^{-s-1}$,  \`a coefficients dans l'anneau $M_{n}(\complex)$ des matrices constantes de dimension $n\times n$,  convergente dans un voisinage de l'infini. Le syst\`eme diff\'erentiel de matrice $A(x)$  est  fuchsien. Il est en outre non r\'esonnant si deux valeurs propres distinctes de la matrice $A_0$ ont une diff\'erence non enti\`ere. Dans ce cas, il existe une (unique) transformation de jauge   tangente \`a l'identit\'e donn\'ee par une s\'erie : 
$$T(x)=I+\sum_{s\geq 1} T_sx^{-s}$$ convergente  au voisinage de l'infini, telle que 
$B(x)=A_0x^{-1}$.
\item Supposons que $\displaystyle A(x)$ est la somme d'une s\'erie $\sum_{s\geq 0}A_s x^{-s}$,  \`a coefficients dans l'anneau $M_{n}(\complex)$ des matrices constantes de dimension $n\times n$, convergente dans un voisinage de l'infini, avec $A_0$ inversible. Le syst\`eme aux $q$-diff\'erences de matrice $A(x)$ est fuchsien. Il est non r\'esonnant si deux valeurs propres distinctes de $A_0$ ne sont pas congrues modulo $q^\rinteger$ et dans ce cas une (unique) transformation de jauge tangente \`a l'identit\'e le transforme en le syst\`eme de matrice constante $A_0$, (voir \cite{sauloy} p. 1033 o\`u le probl\`eme est trait\'e au voisinage de $0$; on s'y ramène par le changement de variable $x \leftarrow 1/x$).
\item Afin d'obtenir des r\'esultats analogues pour les syt\`emes aux diff\'erences il convient de remplacer les s\'eries de puissances par les \og s\'eries de factorielles\fg{}. %Il convient aussi de supposer  le syst\`eme \'ecrit sous la forme 
%$$(x-1)[Y(x)-Y(x-1)]=A(x)Y(x).$$ 
On suppose
$\dis A(x)=I + \sum_{s\geq 0}A_s\frac1{x(x+1)\cdots (x+s)}$. 
Le domaine naturel de convergence d'une telle s\'erie est un demi-plan vertical $\Re x\gg 0$ que l'on peut consid\'erer comme un voisinage de $+\infty$. Comme dans le cas diff\'erentiel, le syst\`eme est dit non r\'esonnant si deux valeurs propres distinctes de $A_0$ ne sont pas congrues modulo $\rinteger$. Dans ce cas il existe une (unique) transformation de jauge tangente \`a l'identit\'e donn\'ee par une s\'erie de factorielles convergente $\displaystyle T(x)=I+\sum_{s\geq 0}T_s\frac1{x(x+1)\cdots (x+s)}$ telle que $B(x)=I-\frac{A_0}{x-1}$. Remarquons que le système aux différences défini par $B(x)$ peut s'écrire : $$(x-1)[Y(x)-Y(x-1)]=A_0 Y(x).$$ 
\end{enumerate}
Dans chaque cas, r\'esoudre le syt\`eme revient alors \`a r\'esoudre un syst\`eme \`a coefficients constants, ce qui se fait en réduisant cette matrice sous forme normale de Jordan puis en introduisant une famille de \og caract\`eres\fg{} et de \og logarithmes\fg{} adapt\'es \`a chaque cas. 

Par ailleurs ces trois familles de syst\`emes sont  reli\'ees entre elles par des propri\'et\'es de \og confluence\fg{}. En effet en notant $\sigma_q$ (\resp $\tau_h$) l'op\'erateur agissant sur une fonction $f$  par $\sigma_qf(x)=f(qx)$ (\resp $\tau_hf(x)=f(x+h)$), on a : 
$$\lim_{q\to 1}\frac{\sigma_q-id}{(q-1)x}=\frac d{dx} $$
et  :
$$\lim_{h\to 0}\frac{\tau_h-id}{h}=\frac d{dx}. $$
L'\'etude de ces confluences, en liaison avec la th\'eorie de Galois, a \'et\'e men\'ee r\'ecemment par J. Sauloy (\cite{sauloy})  pour la premi\`ere famille et par le second auteur (\cite{roques}) pour la seconde. Ces deux auteurs montrent comment r\'ealiser les deux \'etapes (transformation de jauge pour se ramener au cas constant et d\'etermination d'un syst\`eme fondamental de solutions dans le cas constant) d'une fa\c con suffisamment canonique pour qu'un passage \`a la limite soit possible. Dans le premier cas le param\`etre complexe $q$ est de module strictement plus grand que $1$ et tend vers $1$ sur une spirale logarithmique : $q=q_0^\varepsilon$ o\`u $|q_0|>1$ et $\varepsilon $ est un r\'eel positif que l'on fait tendre vers $0$. Dans le second cas $h$ est r\'eel positif et tend vers $0$.\\
Un autre type de confluence a \'et\'e consid\'er\'e par le premier auteur dans \cite{ duv1} :  on part de la famille d'op\'erateurs  d\'efinis par $\sigma_{q,1}f(x)=f(qx+1)$, et on \'etudie la convergence
$\lim_{q\to 1}\sigma_{q,1}=\tau_1$ not\'ee $\tau$. Il s'agit cette fois d'un ph\'enom\`ene de confluence d'op\'erateurs associ\'es \`a une famille d'homographies \`a deux points fixes confluant vers une homographie \`a un point fixe. Rappelons  le fait bien connu que tout syst\`eme $Y(\psi(x))=A(x)Y(x)$ o\`u $\psi$ est une homographie de la sph\`ere de Riemann distincte de l'identit\'e, se ram\`ene par changement de variable homographique soit \`a un syst\`eme aux $q$-diff\'erences (si $\psi$ a deux points fixes) soit \`a un syst\`eme  aux diff\'erences (si $\psi$ a un point fixe). Ainsi le syst\`eme $Y(qx+1)=A_q(x)Y(x)$ devient un syst\`eme aux $q$-diff\'erences si on effectue le changement de variable $x\mapsto (q-1)x+1$. Au niveau de l'op\'erateur, la confluence pr\'ec\'edente est seulement une sp\'ecialisation en $q=1$, mais cette sp\'ecialisation n'est pas possible  dans les solutions obtenues par   le changement de variable indiqu\'e. Dans l'article cit\'e l'\'etude de cette confluence est faite en remarquant qu'il est possible d'exprimer toute s\'erie convergente \`a l'infini en une s\'erie  de la forme 
$\displaystyle a_0+\sum_{s\geq 1}\frac{a_s}{(x;q)_s} $ o\`u $(x;q)_s=\prod_{j=0}^{s-1}(1-q^jx)$ (voir \cite{duv2} p. 344).  En supposant que $q$ est r\'eel et tend vers $1$ par valeurs supérieures, on montre dans \cite{duv1} qu'apr\`es une telle \'ecriture suivie du changement de variable indiqu\'e, les solutions obtenues tendent vers celles du syst\`eme aux diff\'erences limite.

Dans le pr\'esent travail nous g\'en\'eralisons et unifions ces trois articles en consid\'erant une famille de syst\`emes associ\'es \`a la famille d'op\'erateurs $\sigma_{q,h}f(x)=f(qx+h)$. Nous montrons qu'il est possible, sous certaines hypoth\`eses,  
de d\'efinir une  solution canonique pour le syst\`eme :
$$Y(qx+h)=A^{(q,h)}(x)Y(x).$$  Celle-ci est obtenue en r\'esolvant chaque syst\`eme de la famille en les deux \'etapes indiqu\'ees : par une transformation de jauge tangente \`a l'identit\'e d\'efinie par une s\'erie appropriée (séries de $(q,h)$-factorielles) convergente dans un domaine convenable, on se ram\`ene au cas d'un syst\`eme \`a coefficients constants que l'on r\'esout  \`a l'aide d'une famille adapt\'ee de fonctions. De plus ces deux \'etapes sont men\'ees en contr\^olant  la d\'ependance en $(q,h)$. Cette \'etude aboutit au r\'esultat \'enonc\'e dans le Th\'eor\`eme \ref{T1}, section \ref{les deux theos}. Ensuite  on \'etudie les trois types de confluence : $(q,h)\to (q_0,0)$ ($q_0>1$), $(q,h)\to (1,h_0)$ ($h_0>0$) et $(q,h)\to (1,0)$. Les conclusions obtenues sont l'objet du Th\'eor\`eme \ref{T2}, section \ref{les deux theos}.  

Indiquons le plan de l'article. Apr\`es avoir donn\'e les d\'efinitions et notations utiles, nous \'enon\c cons les deux th\'eor\`emes principaux de notre article. Un premier paragraphe \'etudie les s\'eries de $(q,h)$-fractorielles. Chacun des deux paragraphes suivants est consacr\'e \`a la preuve de l'un des th\'eor\`emes.

\subsection{Notations}

 Dans tout ce travail $q$ est un \emph{r\'eel sup\'erieur ou \'egal \`a $1$}  dont l'inverse est not\'e $p$ et $h$ est un \emph{r\'eel positif ou nul}.

Rappelons, pour $q\neq 1$ et $\la\in\complex$, la notation de Jackson : $$[\la]_q=\frac{q^{\la}-1}{q-1}$$ que l'on étend au cas $q=1$ par $[\la]_1=\la$. On pose : $$ \dis\puiss{h}{q}{\la}=h[\la]_q.$$
 On note ${\cal Q}= \mathop]1,+\infty\mathop[\times\mathop]0,+\infty\mathop[$, ${\cal Q}^+=\mathop[1,+\infty\mathop[\times\mathop[0,+\infty\mathop[$ et ${\cal Q}^*={\cal Q}^+ \setminus\{(1,0)\}$.\\
Pour $(q,h)\in {\cal Q}^+$ et $s\in \integer^*$, nous introduisons la fraction rationnelle suivante :
$$
x^{-\puiss{h}{q}{s}}=\frac{1}{x(qx+\puiss{h}{q}{1})\cdots
(q^{s-1}x+\puiss{h}{q}{s-1})}.$$
On \'ecrit aussi $x^{-\puiss{h}{q}{0}}=1$ et on abr\`ege $x^{-\puiss{1}{q}{s}}$ en $x^{-[s]_q}$ et $x^{-\puiss{1}{1}{s}}$ en $x^{-[s]}$. Remarquons que $x^{-\puiss{0}{q}{s}}=q^{-\frac{s(s-1)}2}x^{-s}$.

Afin de faire ressortir les analogies mentionnées dans l'introduction, si $(q,h)\neq (1,0)$,  il sera utile de remplacer l'op\'erateur $\sigma_{q,h}$ par la pseudo-d\'erivation :
$$\Delta_{q,h}=\frac{\sigma_{q,h}^{-1}-id}{\sigma_{q,h}^{-1}(x)-x}=\frac{\sigma_{p,-ph}-id}{(p-1)x-ph}.$$
Cet op\'erateur v\'erifie la \og formule de Leibniz\fg{} suivante :
\begin{eqnarray*}\Delta_{q,h}(fg)(x)&=&\sigma_{p,-ph}f(x)\Delta_{q,h}g(x)+g(x)\Delta_{q,h}f(x)\\
&=&\sigma_{p,-ph}g(x)\Delta_{q,h}f(x)+f(x)\Delta_{q,h}g(x).\end{eqnarray*}
L'analogue de l'op\'erateur \og fuchsien \fg{} $x\frac d{dx}$ est alors l'op\'erateur :
$$\opmixte{p}{h}=p(x-h)\Delta_{q,h}.$$
 On convient naturellement que $\opmixte{1}{0}=x\frac d{dx}$.
 Le r\'esultat suivant est une g\'en\'eralisation imm\'ediate de la Proposition 1.1 de \cite{duv2} .
 
 \begin{lem}
 La fraction rationnelle $x^{-\puiss{h}{q}{s}}$ est un vecteur propre de $\opmixte{p}{h}$ pour la valeur propre $- [s]_q$.
\end{lem}

Nous envisageons dans ce travail des familles de syst\`emes aux $(q,h)$-différences :
\begin{equation} \label{syst}
\opmixte{p}{h}Y(x)=\deform{q}{h}{A}(x)Y(x), \ \ (q,h)\in \Lambda \subset {\cal Q}^+
\end{equation}  o\`u $\deform{q}{h}{A}(x)$ est  une
 série  de la forme :
$$\deform{q}{h}{A}(x) =\sum_{s=0}^{+\infty}A_{s}^{(q,h)}x^{-\puiss{h}{q}{s}}$$
 avec $A_{s}^{(q,h)} \in M_{n}(\complex)$. De telles s\'eries sont appel\'ees s\'eries de $(q,h)$-factorielles formelles à coefficients dans $M_n(\complex)$.
 
 Des conditions suffisantes de convergence de ces s\'eries sont donn\'ees au paragraphe 2 et rendront naturelles les d\'efinitions qui suivent.

  Pour $\lambda, x\in\complex$, la signification habituelle de l'expression   $|x-\frac h{1-q}|>|\lambda-\frac h{1-q}|$ lorsque $(q,h)\in{\cal Q}$ est  prolong\'ee \`a ${\cal Q}^+$ en convenant qu'elle signifie
$\Re(x)>\Re(\lambda)$ si $q=1, \,h>0$ et  $|x|>|\la|$ si $h=0$. Pour $\lambda\in\complex$, on note :
$$V^{(q,h)}(\la)=\left\{ x \in \complex \ | \ \left|x -\frac h{1-q}\right|> \left|\la -\frac h{1-q}\right| \right\}.$$ 
%En remarquant que si $\la_1$ et $\la_2$ sont r\'eels positifs, $\la_2\geq \la_1$ implique $V^{(q,h)}(\la_2)\subset V^{(q,h)}(\la_1)$, on voit 
Notons que la famille $\left(V^{(q,h)}(\la)\right)_{\la\in\real^+}$ constitue un syst\`eme fondamental de voisinages de $\infty$ sauf si $h\neq 0$ et $q=1$ o\`u il s'agit de voisinages de $+\infty$.

On suppose  choisies une norme not\'ee $\|\cdot\|$ sur $M_{n}(\complex)$ et une norme $|\cdot|$ sur $\complex^n$ telles que si $A\in M_{n}(\complex)$ et $X\in \complex^n$, on ait $|AX|\leq \|A\|\,|X|$. Par commodit\'e on supposera aussi que $\|I\|=1$, si $I$ d\'esigne la matrice identit\'e.

La définition suivante fixe le cadre technique dans lequel nous nous placerons.

\begin{defi}
Soient $C$ et $\la$ deux r\'eels positifs et $\Lambda\subset {\cal Q}^+$. Une famille de matrices $\left(\deform{q}{h}{A}(x)\right)_{(q,h)\in\Lambda}$ 
est dite  fuchsienne $(C,\lambda)$ sur $\Lambda$ si, pour tout $(q,h)\in \Lambda$, la
matrice  $\deform{q}{h}{A}(x)$ est une s\'erie de $(q,h)$-factorielles \`a coefficients dans $M_n(\complex)$ :
$$\deform{q}{h}{A}(x)=\sum_{s\geq 0}A_s^{(q,h)}x^{-\puiss{h}{q}{s}}$$ et si, pour tout $s\geq 1$ : $$\|A_s^{(q,h)}\|\leq Cq^{s-1}\left(\la^{-\puiss{h}{q}{s-1}}\right)^{-1}.$$
Une famille de syst\`emes {\em (\ref{syst})} est dite   fuchsienne $(C,\lambda)$ sur $\Lambda$ si : \begin{enumerate}
\item 
pour tout $(q,h)\in \Lambda$,   la matrice $I+(1-q)\deform{q}{h}{A_{0}}$ est inversible,
\item  la famille de matrices $\left(\deform{q}{h}{A}(x)\right)_{(q,h)\in\Lambda}$  est fuchsienne $(C,\lambda)$ sur $\Lambda$.
\end{enumerate}
 Elle
 est dite {\em non
résonnante} si, pour toute valeur propre $c$ de $\deform{q}{h}{A_{0}}$, $q^s c-[s]_q$ n'est valeur propre de $\deform{q}{h}{A_{0}}$ pour aucun entier $s\neq 0$.
\end{defi}

En particulier, chaque système d'une telle famille fuchsienne est fuchsien dans les sens classiques : lorsque $h=0$, la d\'efinition co\"\i ncide avec celle de syst\`eme fuchsien en $\infty$ au sens habituel et lorsque $q=h=1$ c'est la d\'efinition d'un syst\`eme de premi\`ere esp\`ece en $+\infty$ de \cite{harris}.

Notons que la condition d'inversibilit\'e de la matrice $I+(1-q)\deform{q}{h}{A_{0}}$ assure que la matrice du syst\`eme (\ref{syst}) \'ecrit sous la forme $Y(qx+h)=B(x)Y(x)$ est inversible, en accord avec le fait que $\sigma_{q,h}$ est un automorphisme. Elle est \'evidemment r\'ealis\'ee si $q=1$.

Remarquons que la condition de non-r\'esonnance ne fait pas intervenir la valeur de $h$. 
  Lorsque $q=1$, c'est la condition classique : deux valeurs propres distinctes ne diff\'erent pas d'un entier non nul.  Lorsque $q>1$, en remarquant que pour $\tilde c\in\complex$, $q^s[\tilde c]_q+[s]_q=[s+\tilde c]_q$, on voit qu'en notant les valeurs propres de $\deform{q}{h}{A_{0}}$ sous la forme $c=-[-\tilde c]_q$, ce qui est licite vue la condition d'inversibilit\'e indiqu\'ee, la non-r\'esonnance \'equivaut \`a ce que  $\tilde c$ et $\tilde d$ ne sont pas congrus modulo $\rinteger$ lorsque $\tilde c\neq \tilde d$.  On retrouve bien  la condition de non-r\'esonnance rappel\'ee plus haut.

Nous terminons cette section en introduisant la famille de fonctions que nous utiliserons pour r\'esoudre les systèmes à matrice constante. 

Pour tout  $x\in\complex$, si $s\in\integer^*$, on pose $(x;p)_s=\prod_{j=0}^{s-1}(1-p^jx)$ et on  d\'efinit  $$\dis (x;p)_\infty=\lim_{s\to \infty}(x;p)_s.$$ Nous utiliserons aussi la fonction theta de Jacobi suivante dont les propriét\'es de base sont rappelées en \ref{car et log}:
\begin{equation*}
  \Theta_{q}(x)=\sum_{n\in
  \rinteger}(-1)^{n}q^{-\frac{n(n-1)}{2}}x^{n}.
\end{equation*}
Pour $\alpha\in \complex^*$, on note $ \Theta_{q,\alpha}(x)= \Theta_{q}(\alpha^{-1}x)$.

Nous noterons $\Gamma$ la fonction Gamma d'Euler.
 \begin{defi} Pour $(q,h)\in \mathcal{Q}^+$ et $c\in \complex$, on note $ \deform{q}{h}{e_{c}}(x)$ la fonction :  
\begin{equation*}
  \deform{q}{h}{e_{c}}(x)=\left\{\begin{array}{lcl}\dis \left(\frac{h}{1-p}\right)^{c}
  \frac{((1-p)\frac{x}{h}+p;p)_{\infty}}{(p^{c}((1-p)\frac{x}{h}+p);p)_{\infty}}&\mbox{ si } & q>1, h>0\\\dis
  h^{c}\frac{\Gamma(1+c-\frac{x}{h})}{\Gamma(1-\frac{x}{h})}&\mbox{ si } & q=1, h>0\\\dis
  \frac{\Theta_{q}(x)}{\Theta_{q,q^c}(x)}&\mbox{ si } & q>1, h=0\\\dis
  (-x)^c &\mbox{ si } & (q,h)=(1,0)\end{array}\right.
\end{equation*}
Dans le dernier cas, on suppose fix\'ee une d\'etermination du logarithme et on abr\`ege $\deform{1}{0}{e_{c}}$ en $e_c$.\\Pour tout  entier $k\geq 0$, on pose :
$$\deform{q}{h}{\ell_{c,k}}(x)=\frac1{k!}\frac{\partial^k}{\partial c^k} \deform{q}{h}{e_{c}}(x).$$
\end{defi}

Soient $r\geq 1$ un entier et $E_{i,j}$ la matrice carr\'ee \'el\'ementaire de dimension $r$ dont tous les \'el\'ements sont nuls sauf celui situ\'e sur la ligne $i$ et la colonne $j$ qui vaut $1$. On note $\dis N_r=\sum_{i=1}^{r-1}E_{i,i+1}$.
Pour $c\in\complex$ et $(q,h)\in{\cal Q}^+$, on d\'efinit les matrices $r\times r$ :
$$\deform{q}{h}{{\cal L}_{c,r}}(x)=\sum_{j=0}^{r-1}\deform{q}{h}{\ell_{c,j}}(x)  N_r^j
$$ 
et :
 $$\deform{q}{h}{A_{c,r}}=-[-c]_qI_r+\frac{q^{-c}}{1-q}(q^{-N_r}-I_r)$$ 
 
 Soient $r_1,\cdots, r_m$ des entiers strictement positifs de somme $n$ et $c_1,\cdots, c_m$ des nombres complexes, on note : 
 $$\deform{q}{h}{J_{0}}=diag \,(\deform{q}{h}{A_{c_1,r_1}},\cdots, \deform{q}{h}{A_{c_m,r_m}}).$$
  Soit $P$ une matrice constante inversible et $\deform{q}{h}{A_{0}}=P\deform{q}{h}{J_{0}}P^{-1}$. 
\begin{defi} \label{EJ}Avec les notations pr\'ec\'edentes, on  d\'efinit les matrices $n\times n$ : $${\cal E}_{\deform{q}{h}{J_{0}}}(x)= diag \,(\deform{q}{h}{{\cal L}_{c_1,r_1}}(x),\cdots, \deform{q}{h}{{\cal L}_{c_m,r_m }}(x))$$ et :
$${\cal E}_{\deform{q}{h}{A_{0}}}(x)=P{\cal E}_{\deform{q}{h}{J_{0}}}(x)P^{-1}.$$
\end{defi}

\subsection{Les deux th\'eor\`emes.}\label{les deux theos}

Une  transformation de jauge de matrice $F(x)$ transforme le syst\`eme $\opmixte{p}{h}Y(x)=A(x)Y(x)$ en le syst\`eme  de matrice :
$$A^F(x)=F(px-ph)^{-1}[A(x)F(x)-\opmixte{p}{h}F(x)].$$ Il s'agit de la matrice codant le système obtenu à partir de $\opmixte{p}{h}Y(x)=A(x)Y(x)$ par le changement de fonction inconnue  $X(x)=F(x)Y(x)$.
Nous utiliserons des transformations de jauges développables en séries de $(q,h)$-factorielles tangentes \`a l'identit\'e, c'est-\`a dire de la forme :
$${\deform{q}{h}{F}}(x)=I+\sum_{s\geq 1}F_s^{(q,h)}x^{-\puiss{h}{q}{s}}.$$
\'Enon\c cons le premier de nos deux r\'esultats principaux. 
\begin{theo}\label{T1}
Soit $\Lambda$  une partie born\'ee de ${\cal Q}^+$. Soit une famille fuchsienne  non r\'esonnante  $(C,\lambda)$ sur $\Lambda$ de syst\`emes {\em (\ref{syst})}
 o\`u :
 $$\deform{q}{h}{A}(x)=\sum_{s\geq 0}\deform{q}{h}{A_s}x^{-\puiss{h}{q}{s}}.$$ On suppose que, pour tout $(q,h)\in\Lambda$, il existe $\deform{q}{h}{P}\in Gl_n(\complex)$ telle que : $$\deform{q}{h}{P}\deform{q}{h}{J_{0}}=\deform{q}{h}{A_{0}}\deform{q}{h}{P}$$ o\`u :
$$\deform{q}{h}{J_{0}}=diag \,(\deform{q}{h}{A_{c_1,r_1}},\cdots, \deform{q}{h}{A_{c_m,r_m}})$$ avec  $c_1,\cdots, c_m$  des nombres complexes v\'erifiant $c_i-c_j\not\in \rinteger^*$ et $r_1,\cdots,r_m$  des entiers strictement positifs tels que $r_1+\cdots+r_m=n$. Si les familles $\left(\deform{q}{h}{P}\right)_{(q,h)\in \Lambda}$  et $\left(\left(\deform{q}{h}{P}\right)^{-1}\right)_{(q,h)\in \Lambda}$ sont born\'ees,
alors : \begin{enumerate}
\item pour tout $(q,h)\in\Lambda$, il existe une unique transformation de jauge formelle tangente \`a l'identit\'e $\deform{q}{h}{F}(x)$, d\'efinie par une  série de $(q,h)$-factorielles de terme constant $I$,  telle que $\left(\deform{q}{h}{A}\right)^{\deform{q}{h}{F}}(x)=\deform{q}{h}{A_0}$,
\item il existe $(C',\lambda')$ tels que la famille $\left(\deform{q}{h}{F}(x)\right)_{(q,h)\in\Lambda}$ est fuchsienne $(C',\lambda')$,
\item les deux propri\'et\'es suivantes sont satisfaites :
\begin{enumerate}
\item pour toute matrice $P\in Gl_n(\complex)$, $\deform{q}{h}{F}(x)P{\cal E}_{\deform{q}{h}{J_{0}}}(x)$ est une matrice fondamentale de solutions dans $V^{(q,h)}(\la')$  du syst\`eme {\em (\ref{syst})} si et seulement si $P\deform{q}{h}{J_{0}}=\deform{q}{h}{A_{0}}P$,
\item la matrice :
$$\deform{q}{h}{{\cal X}_{can}}(x)=\deform{q}{h}{F}(x){\cal E}_{\deform{q}{h}{A_{0}}}(x)$$ est, dans $V^{(q,h)}(\la')$, une matrice fondamentale de solutions de {\em (\ref{syst})}, ind\'ependante du choix de la matrice $P$ conjuguant $\deform{q}{h}{A_{0}}$ \`a sa forme normale.
\end{enumerate}
\end{enumerate}
\end{theo}

\begin{defi}
La famille de matrices $\left (\deform{q}{h}{{\cal X}_{can}}(x)\right)_{(q,h)\in\Lambda}$ d\'efinie dans le th\'eor\`eme pr\'ec\'edent est appel\'ee solution canonique de la famille {\em (\ref{syst})}.
\end{defi}

Avant d'énoncer des propriétés de confluence pour les fonctions, décrivons le comportement de leurs domaines de définition. Pour $\lambda>0$ fix\'e, la famille $V^{(q,h)}(\la)$ v\'erifie les propri\'et\'es suivantes de continuit\'e en un point du bord de ${\cal Q}^+$~: \begin{enumerate}
\item Pour $q_0>1$ fix\'e, la famille $V^{(q_0,h)}(\la)$ cro\^\i t lorsque $h\downarrow 0$ et $\dis\bigcup_{h>0}V^{(q_0,h)}(\la)=\{|x|>\la\}=V^{(q_0,0)}(\la)$. Notons que cette limite est ind\'ependante de $q_0>1$ et que le r\'esultat est valable si $(q,h)\to(q_0,0)$ avec $q_0>1$.
\item  Pour $h_0>0$ fix\'e, la famille $V^{(q,h_0)}(\la)$ d\'ecro\^\i t lorsque $q\downarrow 1$ et
$\dis\bigcap_{q>1}V^{(q,h_0)}(\la)=\{\Re x>\la\}=V^{(1,h_0)}(\la)$. Cette fois encore le r\'esultat persiste si $(q,h)\to (1,h_0)$ avec $h_0>0$.
\item  La famille $V^{(q,h)}(\la)$ d\'ecro\^\i t avec $\frac h{q-1}$. Chaque  $V^{(q,h)}(\la)$ contient  le demi-plan $\Re x>\la$. Si $\Lambda'\subset \Lambda$, l'intersection des  $V^{(q,h)}(\la)$ pour $(q,h)\in\Lambda'$  est un voisinage de $\infty$ si et seulement si la famille $\left(\frac h{q-1}\right)_{(q,h)\in\Lambda'}$ est born\'ee. Lorsqu'on fait l'hypoth\`ese  suivante : \begin{center}({\bf H})\hspace{0.5cm} $(q,h)\to (1,0)$ de sorte que $\dis \frac h{q-1}\to \ell$ \end{center}o\`u $\ell\in\mathop[0,+\infty\mathop]$, alors $V^{(q,h)}(\la)$ a pour limite $\{|x+\ell|>\la+\ell\}$ si $\ell$ est fini et $\{\Re x>\la\}$ si $\ell=+\infty$. On note $V_\ell(\la)$ ce domaine limite.
\end{enumerate}

\bigskip
 
Nous \'etudierons la confluence des solutions (i.e. la continuité en les paramètres $(q,h)$) dans les trois cas suivants : \begin{description}
\item[C1] $(q,h)$ tend vers $(q_0,0)$ avec $q_0>1$,
\item[C2] $(q,h)$ tend vers $(1,h_0)$ avec $h_0>0$,
\item[C3] $(q,h)$ tend vers $(1,0)$ et l'hypoth\`ese {\bf H} est satisfaite.
\end{description}

\'Enon\c cons le deuxi\`eme r\'esultat principal de cet article.
\begin{theo}\label{T2} Soit une famille {\em(\ref{syst})} v\'erifiant les hypoth\`eses du th\'eor\`eme \ref{T1} dont on reprend les notations. Soit $(q_0,h_0)\in\overline{\Lambda} \cap \partial {\cal Q}^+ $. On suppose que : \begin{enumerate}
\item pour tout $s\geq 0$, $\deform{q}{h}{A}_s$ a une limite $A_s$ lorsque $(q,h)\in\Lambda$ tend vers $(q_0,h_0)$ de telle sorte que l'on soit dans l'une des situations {\em {\bf C1}, {\bf C2}} ou {\em {\bf C3}} et, dans ce dernier cas, on suppose  $\ell \in \real^{+*}$,
\item il existe des matrices de conjugaison $\deform{q}{h}{P}$ telles que, quand $(q,h)\to(q_0,h_0)$, $\deform{q}{h}{P}$   converge vers une matrice conjuguant $\deform{q_0}{h_0}{A_0}$ \`a $\deform{q_0}{h_0}{J_0}$.
\end{enumerate}
  Alors : \begin{enumerate}
\item la s\'erie de $(q_0,h_0)$-factorielles~: $$A(x)=\sum_{s\geq 0}A_sx^{-\puiss{h_0}{q_0}{s}}$$ converge \begin{itemize}
\item dans $V^{(q_0,h_0)}(\la)$ dans les cas {\em {\bf C1}} et {\em {\bf C2}},
\item  dans $V_\ell(\la)$ dans le cas {\em {\bf C3}},
\end{itemize}
\item si le système aux $(q_0,h_0)$-diff\'erences défini par $A(x)$ est fuchsien  et non résonnant, alors
\begin{description}
\item[C1] en posant   $h_n=(1-p_0)p_0^{n}$, la suite $\left( \deform{q_0}{h_n}{{\cal X}_{can}}(x)\right)_{n\geq 0}$   converge  vers $\deform{q_0}{0}{{\cal X}_{can}}(x)$, uniform\'ement sur tout compact de $\{|x|>\la'\} \setminus (c+h_0\integer^*)$,
\item[C2] la famille $\deform{q}{h}{{\cal X}_{can}}(x)$   tend vers $\deform{1}{h_0}{{\cal X}_{can}}(x)$, uniform\'ement sur tout compact de $\{\Re x>\la'\} \setminus (c+h_0\integer^*)$,
\item[C3] la famille $\deform{q}{h}{{\cal X}_{can}}(x)$   tend vers $\deform{1}{0}{{\cal X}_{can}}(x)$, uniform\'ement sur tout compact de $V_\ell(\la') \setminus \real^{+} $.
\end{description}

\end{enumerate}

\end{theo}

Notons que dans le cas {\bf C1} la famille des caractères n'admet pas de limite lorsque $h$ tend vers $0$. En revanche, le théorème \ref{Convc}, section \ref{convcl} assure que cette famille est normale et donne  l'ensemble de ses valeurs d'adhérence. 

\section{Quelques propri\'et\'es des séries de $(q,h)$-factorielles}
Nous donnons ici quelques propriétés des séries de $(q,h)$-factorielles \textit{i.e.} des séries de la forme $\sum_{s\geq 0}A_sx^{-\puiss{h}{q}{s}}$.

\subsection{Conditions suffisantes de convergence}

On a d\'ej\`a remarqu\'e que $x^{-\puiss{0}{q}{s}}=q^{-\frac{s(s-1)}2}x^{-s}$. Les s\'eries de $(q,0)$-factorielles sont donc les s\'eries de puissances en $x^{-1}$. Lorsqu'il n'est pas vide, leur domaine de convergence absolue est de la forme $V^{(q,0)}(\lambda)$ pour un r\'eel positif $\lambda$.

Les s\'eries de $ (1,1)$-factorielles sont les  s\'eries de factorielles. Leur domaine de convergence absolue, s'il est non vide,  est de la forme $V^{(1,1)}(\lambda)$ pour un r\'eel positif $\lambda$ (voir \cite{harris} par exemple).

Enfin, les s\'eries de $(q,1)$-factorielles sont \'etudi\'ees dans \cite{duv2} (p. 351) o\`u elles sont appel\'ees \og s\'eries de factorielles mixtes\fg{}. Toujours sous r\'eserve de converger en au moins un point, elles convergent absolument dans un domaine $V^{(q,1)}(\lambda)$ pour un r\'eel positif $\lambda$.

Pour tout entier $s\geq 1$ et tout $h>0$ :
$$ x^{-\puiss{h}{q}{s}}=h^{-s}\left(\frac xh\right)^{-[s]_q}.$$
Lorsque $h\neq 0$, le changement de variable $x\mapsto \frac xh$ transforme donc une s\'erie de $(q,h)$-factorielles en une s\'erie de $(q,1)$-factorielles. 
On en d\'eduit qu'une s\'erie de $(q,h)$-factorielles qui converge en un point, converge absolument dans un domaine $V^{(q,h)}(\lambda)$ pour un r\'eel positif $\lambda$. 

Par ailleurs, il r\'esulte de \cite{harris} pour les s\'eries de factorielles et de \cite{duv2} (Proposition 2.8) pour les s\'eries de $(q,1)$-factorielles que l'on peut \'enoncer le r\'esultat g\'en\'eral suivant.
\begin{prop}\label{prop1}
Soit $\left(\deform{q}{h}{A}(x)\right)_{(q,h)\in\Lambda}$ une famille
  fuchsienne $(C,\lambda)$ sur $\Lambda\subset {\cal Q}^+$. Pour tout $(q,h)\in\Lambda$ la s\'erie de $(q,h)$-factorielles $\deform{q}{h}{A}(x)$ converge absolument dans $V^{(q,h)}(\lambda)$ et sa somme est une fonction holomorphe.
\end{prop}

 Une fonction $a$ à valeurs complexes, définie et holomorphe sur $V^{(q,h)}(\la)$, est dite d\'e\-ve\-lop\-pable en s\'erie de $(q,h)$-factorielles s'il existe une suite de nombres complexes $(a_s)_{s\geq 0}$ et un r\'eel $\la'>0$ tels que, pour tout $x \in V^{(q,h)}(\la')$~:
$$a(x)=\sum_{s=0}^{+\infty}a_sx^{-\puiss{h}{q}{s}}.$$
 Lorsqu'il existe, un tel développement est unique. L'ensemble des fonctions développables en séries de $(q,h)$-factorielles est noté $\qhfact{q,h}$.

Notons qu'il s'agit d'un anneau commmutatif unitaire et intègre pour l'addition et la multiplication usuelles des fonctions.
De plus, l'application: $$a(x)=\sum_{s=0}^{+\infty}a_sx^{-\puiss{h}{q}{s}} \in \qhfact{q,h} \mapsto \widehat{a}(x)=\sum_{s=0}^{+\infty}a_sx^{-\puiss{h}{q}{s}} \in \complex[[x^{-1}]]$$ est un monomorphisme d'anneaux mais pas un épimorphisme.

L'exemple de  base est   la fonction $\displaystyle \frac 1{x-\la}$ ($\lambda \in \complex$) qui admet le d\'eveloppement en s\'erie de $(q,h)$-factorielles convergente dans $V^{(q,h)}(\la)$ :
\begin{equation}\label{ast}
\frac1{x-\la}=\sum_{s=1}^\infty q^{s-1}\left(\la^{-\puiss{h}{q}{s-1}}\right)^{-1}x^{-\puiss{h}{q}{s}}.
\end{equation}

Remarquons que si $\la$ est  r\'eel positif cette s\'erie est \`a coefficients r\'eels positifs, ce qui permet de l'utiliser comme s\'erie majorante.

\subsection{Formule de translation.}
La combinatoire des s\'eries de factorielles (formelles) est rendue possible par l'existence d'une  \og formule de translation\fg{} permettant d'exprimer une s\'erie en $x^{-[s]}$ comme s\'erie en $(x+\rho)^{-[s]}$. Un r\'esultat analogue (formule d'homoth\'etie) existe pour les s\'eries de $(q,1)$-factorielles (\cite{duv2} p. 352) et le changement de variable $x\to \frac xh$ permet de l'obtenir  pour les s\'eries de $(q,h)$-factorielles avec $h>0$. On prouve en effet que,  si $\displaystyle A(x)=A_0+\sum_{s\geq 1}A_sx^{-\puiss{h}{q}{s}}$, alors :
$$ A(px-ph)=A_0+\sum_{s\geq 1}q^s\left(A_s+[s-1]!\sum_{k=1}^{s-1}h^{s-k}\frac{A_k}{[k-1]!} \right)x^{-\puiss{h}{q}{s}}$$
o\`u   $[0]!=1$  et $[j]!=j!$ si $q=1$,  $[j]!=\prod_{k=1}^j[k]_q$ si $j\geq 1$ (notation de Jackson).\\
Si $h=0$, cette formule se v\'erifie sans recourir \`a une formule de translation.

\subsection{Produit }

Lorsque $h\neq 0$, la formule de multiplication donnée ci-dessous permet d'exprimer le produit de deux s\'eries de $(q,h)$-factorielles formelles comme une série de m\^eme nature. Ces formules figurent dans \cite{harris} dans le cas $(1,1)$, dans \cite{duv2} (Proposition 2.9) dans le cas $(q,1)$ et la formule ci-dessous en résulte. 

Soient  $\dis A(x)= \sum_{s=0}^{+\infty}A_{s}x^{-\puiss{h}{q}{s}}$
et $\dis
B(x)=\sum_{s=0}^{+\infty}B_{s}x^{-\puiss{h}{q}{s}}$ deux séries formelles de $(q,h)$-factorielles.
Leur produit est la s\'erie formelle :
\begin{equation*}
 (AB)(x)=\sum_{s=0}^{+\infty}C_{s}x^{-\puiss{h}{q}{s}}
\end{equation*}
avec:
\begin{equation*}
C_{0}=A_{0}B_{0}; \
C_{s}=A_{0}B_{s}+A_{s}B_{0}+\sum_{(i,j,k)
\in J_{s}
}h^{k}q^{k+ij}c_{i,j,k}A_{i}B_{j},\ s \geq 1
\end{equation*}
o\`u:\begin{equation*}
J_{s} = \{(i,j,k) \;|\; i,j \geq 1,\ k \geq 0,\ i+j+k=s\}
\end{equation*}et:
\begin{equation*}
c_{i,j,k}=\frac{[i+k-1]![j+k-1]!}{[i-1]![j-1]![k]!}
\end{equation*}
 Notons que, pour tout $i,j$, $c_{i,j,0}=1$ et que cette formule se simplifie  en l'habituelle lorsque $h=0$.

Lorsque les s\'eries $A(x)$ et $B(x)$ convergent, leur produit converge au moins dans le plus petit de leurs domaines de convergence. 

\subsection{Inverse}

Pour les s\'eries de $(q,h)$-factorielles \`a coefficients dans $\complex$, nous ferons usage du r\'esultat suivant (\cite{duv2}).

\begin{prop}\label{inverse}Soient $\la,\mu$ des r\'eels positifs et
$\displaystyle a(x)=1-\sum_{s=1}^{+\infty}a_{s}x^{-\puiss{h}{q}{s}}$
une série de $(q,h)$-factorielles dont les  coefficients sont des nombres complexes v\'erifiant $\displaystyle |a_{s}|\leq  \mu q^{s-1}\left(\lambda^{-\puiss{h}{q}{s-1}}\right)^{-1}$,
alors $b=1/a$ est une série de $(q,h)$-factorielles de la forme
$\displaystyle b(x)=1+\sum_{s=1}^{+\infty}{b_{s}}x^{-\puiss{h}{q}{s}}$
avec \\$\displaystyle |b_{s}|\leq  \mu q^{s-1}\left((\lambda+\mu)^{-\puiss{h}{q}{s-1}}\right)^{-1}$.
\end{prop}

\section{Preuve du Th\'eor\`eme 1}

Ce paragraphe est consacré à la preuve du théorème \ref{T1}. 

Comme indiqué plus haut, pour résoudre un système aux $(q,h)$-différences fuchsien, on se ramène au cas constant par
une transformation de jauge convenable. On \'etudie d'abord un syst\`eme de matrice constante (en $x$) non r\'esonnante, puis on montre l'existence d'une unique transformation de jauge ramenant \`a ce cas. En vue de l'étude ultérieure de la confluence, cette deuxi\`eme \'etape est réalisée en famille (famille fuchsienne). 

\subsection{Cas d'une matrice constante}

Résoudre un système à matrice constante, revient  à résoudre les deux familles d'équations suivantes : 
\begin{equation}\label{C1}
  \opmixte{p}{h} y(x) =-\puiss{}{q}{-c}y(x) \ \ \ \text{(\og équation des caractères\fg{})}
\end{equation}
et :
\begin{equation}\label{L1}\opmixte{p}{h} \deform{q}{h}{\ell_{c,k}}(x)=-\puiss{}{q}{-c}\deform{q}{h}{\ell_{c,k}}(x) +\sum_{j=1}^k\frac{(-1)^j}{j!}\frac{\ln ^jq}{1-q}q^{-c}\deform{q}{h}{\ell_{c,k-j}}(x) \ \ \ \text{(\og équation des logarithmes\fg{})}.\end{equation}

\subsubsection{Caractères et  logarithme} \label{car et log}

Les fonctions $ \deform{q}{h}{e_{c}}(x)$ d\'efinies dans l'introduction sont m\'e\-ro\-mor\-phes sur $\complex$ lorsque  $(q,h)\neq (1,0)$. Rappelons quelques propriétés fondamentales des fonctions spéciales utilisées pour les construire (à savoir : la fonction $\Gamma$, le symbole de Pochhammer et  la fonction $\Theta_q$). La fonction $\Gamma$ est m\'eromorphe sur $\complex $, sans z\'eros,  \`a p\^oles simples aux entiers n\'egatifs et v\'erifie $\Gamma(x+1)=x\Gamma(x)$. La fonction $(x;p)_\infty$ est enti\`ere, s'annule aux points  $q^\integer$ et v\'erifie $(x;p)_\infty=(1-x)(px;p)_\infty$. 
La fonction  $ \Theta_{q}(x)$ est holomorphe sur $\complex^{*}$ et v\'erifie  $ \Theta_{q}(qx)=qx\Theta_{q}(x)$. Rappelons aussi la formule du triple
produit de Jacobi :
\begin{equation*}
  \Theta_{q}(x)=(p;p)_{\infty}(x;p)_{\infty}(px^{-1};p)_{\infty}
\end{equation*}  qui montre que $ \Theta_{q}(x)$  a pour ensemble de z\'eros $q^{\rinteger}$. Ces propri\'et\'es permettent d'\'etablir la proposition qui suit.

\begin{prop}\label{carac} Pour $(q,h)\in \mathcal{Q}^+$ et $c\in \complex$, la fonction $ \deform{q}{h}{e_{c}}(x)$ est solution de l'équation des caractères {\em (\ref{C1})}.
Si $c=-s$ est un entier {\em n\'egatif}, on a :
$$\deform{q}{h}{e_{-s}}(x)=(-1)^sq^{-s}x^{-\puiss{h}{q}{s}}.$$
Si $c\in\integer$, $ \deform{q}{h}{e_{c}}(x)$ est un polyn\^ome de degr\'e $c$. \\Enfin, si $c\not\in\integer$ et si $(q,h)\neq(1,0)$, les p\^oles de $ \deform{q}{h}{e_{c}}(x)$ sont simples et situ\'es aux points suivants :\begin{enumerate}
\item[(i)] $h[c+\integer^*]_q=\{h[c+n]_q\;|\;n\in\integer^*\}$ si $q>1$ et $h>0$,
\item[(ii)] $c+\rinteger$ si $q>1$ et $h=0$,
\item[(iii)] $h(c+\integer^*)$ si $q=1$ et $h>0$.
\end{enumerate}

\end{prop}

Remarquons que ce  choix de caract\`eres se sp\'ecialise en celui de \cite{duv1} lorsque $h=1$, celui de \cite{roques} lorsque $q=1$ et celui de \cite{sauloy} lorsque $h=0$, en remarquant dans ce dernier cas que la valeur propre est ici not\'ee $q^c$. De plus,
\`a constante près, lorsque $q\neq 1$ et $h\neq 0$, cette solution de (\ref{C1})
coïncide avec la solution canonique utilis\'ee dans \cite{sauloy} pour r\'esoudre l'\'equation aux $q$-différences en $t$ (régulière en $0$) obtenue à partir de (\ref{C1}) après le changement de variable
$t=1+(q-1)\frac xh$.

Ces fonctions permettent, comme on le verra plus bas, de résoudre les systèmes à matrice constante semi-simple. Pour traiter le cas g\'en\'eral on s'appuie sur le r\'esultat suivant qui s'obtient par d\'erivations successives  par rapport \`a $c$  de  (\ref{C1}), en convenant que, si $q=1$, $\frac{\ln q}{q-1}q^{-c}=1$ et, pour $j\geq 2$, $\frac{\ln^j q}{q-1}q^{-c}=0$.

\begin{prop} Pour tout $c\in\complex$ et tout entier $k\geq 1$, les fonctions $\deform{q}{h}{\ell_{c,k}}(x)$  v\'erifient  la famille d'\'equations (\ref{L1}).
\end{prop}
 Si $q=1$, la formule se r\'eduit \`a la relation  :
 $$\opmixte{1}{h} \deform{1}{h}{\ell_{c,k}}(x)=c\deform{1}{h}{\ell_{c,k}}(x) + \deform{1}{h}{\ell_{c,k-1}}(x) $$ et la famille de solutions choisie est  celle utilis\'ee dans \cite{roques}.
La sp\'ecialisation $h=1$ correspond \`a l'une des deux familles utilis\'ees dans \cite{duv1} (celle not\'ee $\widetilde L_k$). La sp\'ecialisation $h=0$, $q>1$ donne une famille diff\'erente de  celle utilis\'ee dans \cite{sauloy} et dont le lien avec cette autre famille \og naturelle\fg{} est \'etudi\'e dans \cite{duv3}. Enfin la sp\'ecialisation $(q,h)=(1,0)$ correspond \`a la famille habituelle des logarithmes puisque :
$$\deform{1}{0}{\ell_{c,k}}(x)=\frac{\ln^k(-x)}{k!}\,(-x)^c.$$

 \subsubsection{Solution canonique (cas constant)}
 Rappelons la convention  $\frac{\ln q}{q-1}q^{-c}=1$ si $q=1$. Avec les notations de l'introduction, posons :
 \begin{eqnarray}\label{J1}\tilde N_{r_j}=\frac{q^{-c}}{1-q}(q^{-N_{r_j}}-I_{r_j}).\end{eqnarray}  En remarquant que $\dis q^{-N_r}-I_r=\sum_{j=1}^{r-1}(-1)^i\frac{\ln^j q}{j!} N_r^j$, on d\'eduit des formules (\ref{C1}) et (\ref{L1}) la proposition suivante dont la preuve est laissée au lecteur.
\begin{prop}Pour tout $r\geq 1$, la matrice $r\times r$ :
$$\deform{q}{h}{{\cal L}_{c,r}}(x)=
\left(
\begin{array}{ccccc}
  \deform{q}{h}{\ell_{c,0}}(x)&\deform{q}{h}{\ell_{c,1}}(x)   &\deform{q}{h}{\ell_{c,2}}(x) &\cdots & \deform{q}{h}{\ell_{c,r-1}}(x) \\
 0 &\deform{q}{h}{\ell_{c,0}}(x)   & \deform{q}{h}{\ell_{c,1}}(x) &\cdots &\deform{q}{h}{\ell_{c,r-2}}(x) \\
 \vdots  & \vdots  &  \ddots & &\vdots\\
 0&&\cdots &0& \deform{q}{h}{\ell_{c,0}}(x)
\end{array}
\right)
$$ est un syst\`eme fondamental de solutions du syst\`eme de dimension $r$ :
$$\opmixte{p}{h}Y(x)= \deform{q}{h}{A_{c,r}} Y(x).$$ 
\end{prop}

Reprenons les notations de la d\'efinition \ref{EJ}. 

\begin{lem}\label{L10}
Toute matrice constante qui commute avec $\deform{q}{h}{J_{0}}$ commute avec ${\cal E}_{\deform{q}{h}{J_{0}}}(x)$.
\end{lem}

\begin{proof} Commen\c cons par montrer que les matrices cons\-tan\-tes qui commutent avec $\deform{q}{h}{J_{0}}$ sont exactement celles qui commutent avec la forme de Jordan habituelle qui correspond \`a  des blocs de la forme $-[-c]_qI_r+N_r$ . Si $q=1$, les deux formes normales co\"\i ncident et il n'y a rien \`a prouver. \\Si $q>1$,  partageons une  matrice $P$ qui commute avec $\deform{q}{h}{J_{0}}$  en $P=(P_{j\,k})_{1\leq j,k\leq m}$ o\`u $P_{j\,k}$ est une matrice $r_j\times r_k$.  Elle commute avec $\deform{q}{h}{J_{0}}$ si et seulement si $P_{j\,k}=0$ lorsque $[-c_j]_q\neq [-c_k]_q$ et $\tilde N_{r_j}P_{j\,k}=P_{jk}\tilde N_{r_k}$ lorsque $[-c_j]_q=[-c_k]_q$, d'o\`u aussi, pour tout entier $i\geq 1$, $\dis\tilde N_{r_j}^iP_{j\,k}=P_{j\,k}\tilde N_{r_k}^i$. Dans ce cas, si  $c$ est tel que $q^c=q^{c_j}=q^{c_k}$,   on  d\'eduit de (\ref{J1}) que :
\begin{eqnarray}\label{J2}N_{r_j}=-\frac1{\ln q}\ln(I_{r_j}-(q-1)q^c\tilde N_{r_j})=\frac1{\ln q}\sum_{i=1}^{r_j-1}\frac{(q-1)^iq^{ic}}i\tilde N_{r_j}^i\end{eqnarray} et donc aussi : \begin{eqnarray*}
 N_{r_j}P_{j\,k}&=&\frac1{\ln q}\sum_{i\geq 1}\frac{(q-1)^iq^{ic}}i\tilde N_{r_j}^iP_{j\,k}\\&=&\frac1{\ln q}\sum_{i\geq 1}\frac{(q-1)^iq^{ic}}iP_{j\,k}\tilde N_{r_k}^i=P_{j\,k}N_{r_k}.
\end{eqnarray*}

La matrice ${\cal E}_{\deform{q}{h}{J_{0}}}(x)$ est diagonale par blocs correspondants \`a ceux de $\deform{q}{h}{J_{0}}$, le $j$-i\`eme bloc \'etant :
$$\deform{q}{h}{{\cal L}_{c_j,r_j}}(x) =\deform{q}{h}{e_{c_j}}(x)I_{r_j}+\sum_{i=1}^{r_j-1}\deform{q}{h}{\ell_{c_j,i}}(x)N_{r_j}^i.$$  Elle commute avec $P$ si et seulement si pour tout $k,j$ on a $\deform{q}{h}{{\cal L}_{c_j,r_j}}(x) P_{j\,k}=P_{j\,k}\deform{q}{h}{{\cal L}_{c_k,r_k}}(x)$. Lorsque $[-c_j]_q\neq [-c_k]_q$ cette relation est claire puisque $P_{j\,k}=0$. Lorsque $[-c_j]_q=[-c_k]_q$, elle se d\'eduit ais\'ement des relations $N_{r_j}P_{j\,k}=P_{j\,k}N_{r_k}$. \end{proof}

Soit $\deform{q}{h}{A_{0}}$ la matrice du syst\`eme (\ref{syst}). Puisqu'elle est inversible, on peut 
supposer  que ses valeurs propres sont   \'ecrites sous la forme $-[-c_1]_q,\cdots ,-[-c_n]_q$ (r\'ep\'et\'ees s'il y a plusieurs blocs de Jordan pour une m\^eme valeur propre) et qu'elles v\'erifient la condition de non r\'e\-son\-nan\-ce qui s'\'ecrit maintenant  $c_i-c_j\not\in \rinteger^*$. La relation (\ref{J2}) permet de prendre 
 $\deform{q}{h}{J_{0}}$ pour forme normale de $\deform{q}{h}{A_{0}}$, c'est-\`a-dire de supposer qu'il existe $P\in GL_n(\complex)$ telle que $\deform{q}{h}{A_{0}}=P\deform{q}{h}{J_{0}}P^{-1}$.  La matrice suivante :
 $$P{\cal E}_{\deform{q}{h}{J_{0}}}(x)P^{-1}$$ est clairement un système fondamental de solutions du système défini par $\deform{q}{h}{A_{0}}$. Voyons ce qu'il en est de la dépendance de ce système fondamental par rapport à la matrice de passage $P$.
 \begin{lem} Si $P_1$ et $P_2$ sont deux matrices conjuguant $\deform{q}{h}{A_{0}}$ \`a $\deform{q}{h}{J_{0}}$, alors :
 $$P_1{\cal E}_{\deform{q}{h}{J_{0}}}(x)P_1^{-1}=P_2{\cal E}_{\deform{q}{h}{J_{0}}}(x)P_2^{-1}.$$
 \end{lem}
\begin{proof}
L'hypoth\`ese implique que la matrice $P_1^{-1}P_2$ commute avec $\deform{q}{h}{J_{0}}$ et donc, d'apr\`es le lemme \ref{L10} :
$$P_1^{-1}P_2{\cal E}_{\deform{q}{h}{J_{0}}}(x)={\cal E}_{\deform{q}{h}{J_{0}}}(x)P_1^{-1}P_2 $$ ou encore :$$P_1{\cal E}_{\deform{q}{h}{J_{0}}}(x)P_1^{-1}=P_2{\cal E}_{\deform{q}{h}{J_{0}}}(x)P_2^{-1}.$$
\end{proof}

Cela montre que ${\cal E}_{\deform{q}{h}{A_{0}}}(x)$ est effectivement indépendante des matrices conjuguantes choisies et justifie la définition \ref{EJ} : dans le cas d'un syst\`eme de matrice constante, on appelle {\em solution canonique} la matrice ${\cal E}_{\deform{q}{h}{A_{0}}}(x)$.

\subsection{Transformation de jauge}
  Nous allons maintenant \'etablir l'existence d'une (unique) famille de transformations de jauge tangentes \`a l'identit\'e ramenant la famille de syst\`emes (\ref{syst})  \`a une famille de systèmes à coefficients constants, ce dernier cas ayant été traité dans la section précédente. Cela se r\'ealise en plusieurs \'etapes et le th\'eor\`eme qui suit est un outil essentiel.
\begin{theo}\label{outil}
Soit $\Lambda$  une partie born\'ee de ${\cal Q}^+$. Consid\'erons une famille {\em (\ref{syst})}  fuchsienne non résonnante $(C,\lambda)$ sur $\Lambda$. On note $\deform{q}{h}{C_{0}}$ l'inverse de la matrice $I+(1-q)\deform{q}{h}{A_{0}}$ et on suppose qu'il existe $\al>0$ tel que, pour tout $(q,h)\in \Lambda$, $\|\deform{q}{h}{C_{0}}\|\leq \al$. La non r\'esonnance implique que, pour tout $s\geq 1$, l'op\'erateur d\'efini sur $M_n(\complex)$ par :  $U\mapsto  (\deform{q}{h}{A_{0}}+[s]_qI)U-q^sU\deform{q}{h}{A_{0}}$ est inversible. On note $\deform{q}{h}{\widetilde L_s}  $ son inverse et on suppose qu'il existe $M>0$ tel que pour tout $s\geq 1$ et tout $(q,h)\in\Lambda$, l'op\'erateur $\deform{q}{h}{\widetilde L_s}  $ soit de norme inf\'erieure ou \'egale \`a $M$.

Alors,
pour chaque $(q,h)\in\Lambda$, il existe une unique transformation de jauge
$\deform{q}{h}{F}(x)$ dé\-ve\-lop\-pable en série de $(q,h)$-factorielles,
tangente à l'identité, qui transforme le système
précédent en celui de matrice constante
$\deform{q}{h}{A_{0}}$ autrement dit :
$(\deform{q}{h}{A})^{\deform{q}{h}{F}}=\deform{q}{h}{A_{0}}.$ De
plus, il existe  $C'>0$ et $\lambda'>0$ tels que la famille  $(\deform{q}{h}{F}(x))_{(q,h)\in \Lambda}$ est de type
$(C',\lambda')$. 
\end{theo}
Le cas $h=0$ et $q=1$ est le classique th\'eor\`eme concernant les syst\`emes diff\'erentiels fuchsiens.
Les trois sp\'ecialisations $h=1$ ou $q=1$  (et $h>0$) ou  $h=0$ (et $q>1$) donnent les r\'esultats \'etablis dans \cite{duv1}, \cite{roques} ou \cite{sauloy} et la preuve qui suit reprend la m\^eme d\'emarche.

\subsubsection{Un r\'esultat pr\'eliminaire}

La preuve du th\'eor\`eme \ref{T1} repose sur la proposition suivante.

\begin{prop}\label{preparatif}
Soit une famille {\em (\ref{syst})} 
fuchsienne $(C,\lambda)$ sur $\Lambda$. On suppose que, pour tout $s\geq 1$,
$-[s]_q$ n'est pas valeur propre de la matrice $A_0^{(q,h)}$.  Soit $U_{0}$ un vecteur non nul de $\complex^{n}.$ Pour tout $(q,h)\in \Lambda$, le  syst\`eme  {\em (\ref{syst})} a une solution qui est une s\'erie (formelle) de $(q,h)$-factorielles de terme constant
$U_{0}$  si et seulement si $U_{0}$ appartient au noyau de $A_0^{(q,h)}$. Une telle solution  est alors unique. De plus, si
$\displaystyle b=\sup_{s\geq 1,(q,h)\in\Lambda} \left\|\left(\puiss{}{q}{s}I+A_{0}^{(q,h)}\right)^{-1}\right\|$ est
fini,  la famille form\'ee par ces solutions est de type
$(Cb|U_{0}|,Cb+\lambda).$
\end{prop}

\begin{proof}
Notons $\displaystyle \deform{q}{h}{A}(x)=\sum_{s=0}^{+\infty}
A_{s}^{(q,h)}x^{-\puiss{h}{q}{s}}$ et  cherchons $\deform{q}{h}{Y}(x)$ sous la forme $\displaystyle \deform{q}{h}{Y}(x)=U_0+\sum_{s=1}^{+\infty}
Y_{s}^{(q,h)}x^{-\puiss{h}{q}{s}}.$

\textsl{Partie formelle. } D'apr\`es la formule du produit:
\begin{equation*}
\deform{q}{h}{A}(x)\deform{q}{h}{Y}(x)=\sum_{s=0}^{+\infty}C_{s}^{(q,h)}x^{-\puiss{h}{q}{s}}
\end{equation*}
avec
$C_{0}^{(q,h)}=A_{0}^{(q,h)}U_0 $ et, pour $s\geq 1$:
\begin{equation*}
C_{s}^{(q,h)}=A_{0}^{(q,h)}Y_{s}^{(q,h)}+A_{s}^{(q,h)}U_0+\sum_{(i,j,k)
\in J_{s}
}h^kq^{k+ij}c_{i,j,k}A_{i}^{(q,h)}Y_{j}^{(q,h)}.
\end{equation*}
D'autre part, puisque:
\begin{equation*}
\opmixte{p}{h}\deform{q}{h}{Y}(x)=\sum_{s=1}^{+\infty}
-[s]_q
Y_{s}^{(q,h)}x^{-\puiss{h}{q}{s}},
\end{equation*}
 $\deform{q}{h}{Y}$ est solution de (\ref{syst}) si et seulement si ses coefficients v\'erifient :
 \begin{equation*}
   \begin{cases}
A_{0}^{(q,h)}U_{0}=0  \\
 \displaystyle   -[s]_qY_{s}^{(q,h}=A_{0}^{(q,h)}Y_{s}^{(q,h)}+
A_{s}^{(q,h)}U_0+\sum_{(i,j,k)
\in J_{s}
}h^kq^{k+ij}c_{i,j,k}A_{i}^{(q,h)}Y_{j}^{(q,h)}
  \end{cases}
  \end{equation*}
ou encore:
\begin{equation*}
\begin{cases}
 A_{0}^{(q,h)}U_{0}=0  \\
   \displaystyle -([s]_qI+A_{0}^{(q,h)})Y_{s}^{(q,h)}=
A_{s}^{(q,h)}U_0+\sum_{(i,j,k)
\in J_{s}
}h^kq^{k+ij}c_{i,j,k}A_{i}^{(q,h)}Y_{j}^{(q,h)}.
  \end{cases}
  \end{equation*}
Il est donc nécessaire que $U_{0}$ soit dans
$ker A_{0}^{(q,h)}.$ Comme, par hypothèse, la matrice
$[s]_qI+A_{0}^{(q,h)}$ est inversible, le
système se résout de proche en proche :
\begin{equation*}
Y_{s}^{(q,h)}=-\left([s]_qI+A_{0}^{(q,h)}\right)^{-1}
    \left[A_{s}^{(q,h)}U_0+\sum_{(i,j,k)
\in J_{s}
}h^kq^{k+ij}c_{i,j,k}A_{i}^{(q,h)}Y_{j}^{(q,h)}
\right].
  \end{equation*}
Ceci termine l'aspect formel de la proposition.\\
\textsl{Partie convergente. } \'Etablissons maintenant, sous l'hypoth\`ese indiqu\'ee,
la convergence de cette solution.

Notons : $a_{s}^{(q,h)}=\|A_{s}^{(q,h)}\|$, $u_0=|U_0|$ et, $\forall s \geq 1$,
$y_{s}^{(q,h)}=|Y_{s}^{(q,h)}|$. \\
De l'équation récurrente ci-dessus  et de la d\'efinition de $b$, on déduit, pour $s\geq 1$, les majorations :
\begin{equation*}
    \deform{q}{h}{y_{s}}\leq b\left[\deform{q}{h}{a_{s}}u_0
    +\sum_{(i,j,k)
\in J_{s}
}h^kq^{k+ij}c_{i,j,k}\deform{q}{h}{a_{i}}\deform{q}{h}{y_{j}}
\right].
  \end{equation*}
Introduisons la suite d\'efinie r\'ecursivement par :
\begin{equation*}
  \begin{cases}
   \deform{q}{h}{z_{1}}=b\deform{q}{h}{a_{1}}u_0\\
    \displaystyle \deform{q}{h}{z_{s}}= b\left[\deform{q}{h}{a_{s}}u_0
    +\sum_{(i,j,k)
\in J_{s}
}q^{k+ij}h^kc_{i,j,k}\deform{q}{h}{a_{i}}\deform{q}{h}{z_{j}}
\right]
  \end{cases}
  \end{equation*}
de sorte que, pour tout $s\geq 1$,  $\deform{q}{h}{y_{s}}\leq \deform{q}{h}{z_{s}}$. Notons:
\begin{equation*}
  \deform{q}{h}{z}(x)=\sum_{s=1}^{+\infty}
\deform{q}{h}{z_{s}}x^{-\puiss{h}{q}{s}}
  \end{equation*}
et:
\begin{equation*}
 \deform{q}{h}{a}(x)=\sum_{s=1}^{+\infty}
\deform{q}{h}{a_{s}}x^{-\puiss{h}{q}{s}}.
  \end{equation*}
La relation suivante est une cons\'equence immédiate de la r\'ecurrence d\'efinissant $\deform{q}{h}{z}(x)$:
$$\deform{q}{h}{z}(x)=b(u_0\deform{q}{h}{a}(x)
+\deform{q}{h}{a}(x)\deform{q}{h}{z}(x)).$$ On en d\'eduit :\begin{equation*}
\deform{q}{h}{z}(x)=\frac{bu_0\deform{q}{h}{a}(x)}{1-b\deform{q}{h}{a}(x)}
=-u_0\left(
1-\frac{1}{1-b\deform{q}{h}{a}(x)}\right).
  \end{equation*}
et il ne reste qu'\`a  appliquer la proposition \ref{inverse} pour obtenir :
\begin{equation*}
\deform{q}{h}{z_{s}} \leq Cbu_0
q^{s-1}\left((\lambda+Cb)^{-\puiss{h}{q}{s-1}}\right)^{-1}.
  \end{equation*}
\end{proof}

Lorsque $\deform{q}{h}{A_{0}}=0$, le systèmes est \og régulier\fg{}.
 En appliquant le résultat précédent  \`a la base canonique $(\ve_1,\cdots, \ve_n)$ de $\complex^n$, on obtient le
\begin{coro} \label{corollaire sur resolution dans le cas regulier}
Soit $\opmixte{p}{h}Y(x)=\deform{q}{h}{A}Y(x)$ une famille de systèmes
réguliers, de type  $(C,\lambda)$ sur $\Lambda$. Ces équations admettent
chacune un système fondamental de solutions, tangent à l'identité,  form\'e de séries de $(q,h)$-factorielles.
La famille ainsi d\'efinie  est de type $(C',C+\lambda)$ avec $\dis C'=C\max_{1\leq i\leq n}|\ve_i|$.
\end{coro}

\subsubsection{Preuve du théorème \ref{outil}}
Notons $\deform{q}{h}{F}(x)=I+\sum_{s=1}^{+\infty}
\deform{q}{h}{F_{s}}x^{-\puiss{h}{q}{s}}$ la transformation de jauge cherch\'ee. La condition qu'elle doit v\'erifier s'\'ecrit :
\begin{equation}\label{az}
   \deform{q}{h}{A}(x)\deform{q}{h}{F}(x)-\opmixte{p}{h}
\deform{q}{h}{F}(x)=\deform{q}{h}{F}(px-ph)\deform{q}{h}{A_{0}}
\end{equation}
ou encore :
\begin{equation*} \label{la transfo de jauge}
  \opmixte{p}{h}
\deform{q}{h}{F}(x)\left[I+\frac{(p-1)x-ph}{p(x-h)}\deform{q}{h}{A_{0}}\right]=\left[\deform{q}{h}{A}(x)\deform{q}{h}{F}(x)-\deform{q}{h}{F}(x) \deform{q}{h}{A_{0}}\right] .
\end{equation*}
%La matrice $B_0^{(q,h)}=I+(1-q)\deform{q}{h}{A_{0}}$ est la limite quand $x\to+ \infty$ de la matrice $B(x)$ du syst\`eme (\ref{syst}) suppos\'e \'ecrit sous la forme $$Y(x)=B(x)Y(qx+h).$$ L'inversibilit\'e de $B(x)$ signifie que le syst\`eme est \og de premi\`ere esp\`ece\fg{}. Les hypoth\`eses du th\'eor\`eme assurent que  cette condition est remplie pour tout $(q,h)\in\Lambda$.

D'autre part, si $a\in \complex$ est tel que $b=1+(1-q)a\neq 0$,  la formule (\ref{ast}) permet d'\'ecrire, en posant $\displaystyle d=\frac{a}{b}$ et $\displaystyle \mu=\frac{(1+a)h}b=(1+qd)h$ :
\begin{eqnarray*}
\left[1+\frac{(p-1)x-ph}{p(x-h)}a\right]^{-1}&=&\frac{x-h}{bx-(1+a)h}=
\frac1{b}\left[1+\frac{qah}{bx-(1+a)h}\right]=\\= \frac1{b}\left[1+\frac{qdh}{x-\mu}\right]&=&
 \frac1{b}\left[1+qdh\sum_{s\geq 1}q^{s-1}\left(\mu^{-\puiss{h}{q}{s-1}}\right)^{-1}x^{-\puiss{h}{q}{s}}\right]\\&=&
 \frac1b\left[1+\sum_{s\geq 1}q^sh^s(d^{-[s]_q})^{-1}x^{-\puiss{h}{q}{s}}\right].
\end{eqnarray*}
 Cette derni\`ere \'egalit\'e s'obtient en remarquant que, puisque $\mu=(1+qd)h$, on a $ dh\left(\mu^{-\puiss{h}{q}{s-1}})\right)^{-1} = h^s(d^{-[s]_q})^{-1}$.

Dans ces expressions, qui sont des polyn\^omes en $d$, on peut remplacer le nombre complexe $d$ par la matrice :  \begin{equation*}
   \deform{q}{h}{D_{0}}=\deform{q}{h}{C_{0}}\deform{q}{h}{A_{0}},
\end{equation*} et donc $a$ par $\deform{q}{h}{A_{0}}$ et $1/b$ par $\deform{q}{h}{C_{0}}$. On obtient  le d\'eveloppement suivant :
\begin{equation*}
   \left[I+\frac{(p-1)x-ph}{p(x-h)}\deform{q}{h}{A_{0}}\right]^{-1}=\sum_{s=0}^{+\infty}\deform{q}{h}{C_{s}}x^{-\puiss{h}{q}{s}}
\end{equation*}
avec, pour $s\geq 1$:
\begin{equation*}
   \deform{q}{h}{C_{s}}=h^sq^{s}\deform{q}{h}{C_{0}} \prod_{k=0}^{s-1}\left(q^k\deform{q}{h}{D_{0}}+[k]_qI\right).
\end{equation*}

L'\'equation (\ref{az}) peut se mettre sous la forme :
\begin{equation*}
  \opmixte{p}{h}
\deform{q}{h}{F}(x)=\left[\deform{q}{h}{A}(x)\deform{q}{h}{F(x)}-\deform{q}{h}{F}(x) \deform{q}{h}{A_{0}}\right] \left[I+\frac{(p-1)x-ph}{p(x-h)}\deform{q}{h}{A_{0}}\right]^{-1}
\end{equation*}
et s'interpréter  comme un
système aux $(q,h)$-différences de dimension $n^2$ auquel on peut appliquer la proposition  \ref{preparatif}. En effet :
 \begin{equation*}
 \opmixte{p}{h}
\deform{q}{h}{F}(x)=\sum_{s=0}^{+\infty}\deform{q}{h}{L}_{s}(\deform{q}{h}{F}(x)) x^{-\puiss{h}{q}{s}}
\end{equation*}
où $\deform{q}{h}{L}_{s}$ est l'opérateur linéaire agissant sur $M_n(\complex)$  par :
\begin{equation*}
  \deform{q}{h}{L}_{0}(U)=\left[\deform{q}{h}{A_{0}}U-U\deform{q}{h}{A_{0}}\right]\deform{q}{h}{C_{0}},
\end{equation*}
et, pour $s \geq 1$:
\begin{eqnarray*}
  \deform{q}{h}{L}_{s}(U)&=&\deform{q}{h}{A_{s}}U\deform{q}{h}{C_{0}}+(\deform{q}{h}{A_{0}}U-U\deform{q}{h}{A_{0}})\deform{q}{h}{C_{s}}+\\&&+
\sum_{(i,j,k) \in J_{s}}h^kq^{k+ij}
c_{i,j,k}\deform{q}{h}{A_{i}}U\deform{q}{h}{C_{j}}.
\end{eqnarray*}
L'opérateur linéaire $\deform{q}{h}{L}_{0}$  admet $0$ pour
valeur propre et $I$ est un vecteur propre associé. Pour tout $s\geq 1$, $-[s]_q$ n'est pas valeur propre de $\deform{q}{h}{L}_{0}$ : cela résulte de l'hypoth\`ese de non r\'esonnance puisque l'\'egalit\'e $ \deform{q}{h}{L}_{0}(U)=-[s]_qU$ \'equivaut \`a : $$\deform{q}{h}{A_{0}}U-U\deform{q}{h}{A_{0}}=-[s]_qU(I+(1-q)\deform{q}{h}{A_{0}})=-[s]_qU+(q^s-1)U\deform{q}{h}{A_{0}}$$ et donc aussi \`a :
$$(\deform{q}{h}{A_{0}}+[s]_qI)U-q^sU\deform{q}{h}{A_{0}}=0.$$
De plus l'hypoth\`ese faite permet de majorer uniform\'ement sur $\Lambda$ la norme de l'op\'erateur $(\deform{q}{h}{L}_0+[s]_qId)^{-1}$ qui s'exprime en fonction de l'op\'erateur $\deform{q}{h}{\widetilde L_s}$ et de la matrice $I+(1-q)\deform{q}{h}{A_{0}}$ dont la norme est born\'ee par hypoth\`ese. Ensuite,
puisque $\|\deform{q}{h}{D_{0}}\|\leq C\al$, on d\'eduit de la formule donnant $\deform{q}{h}{C_{s}}$ la majoration:
\begin{equation*}
\|\deform{q}{h}{C_{s}}\| \leq \al h^s q^{s}\prod_{k=0}^{s-1}(q^kC\al+[k]_q)=\al h^s q^s\left((C\al)^{-[s]_q}\right)^{-1}.
\end{equation*}
D'o\`u l'estimation suivante pour la norme de $\deform{q}{h}{L}_s$ :
\begin{eqnarray*}
\|\deform{q}{h}{L} _{s} \| &\leq&
C\al q^{s-1}\left(\lambda^{-\puiss{h}{q}{s-1}}\right)^{-1}+2C\al h^sq^{s}\left((C\alpha)^{-[s]_q}\right)^{-1}+\\
&& +C\sum_{(i,j,k) \in J_{s}}h^{k+j}q^{ij+s-1}
c_{i,j,k}\left(\lambda^{-\puiss{h}{q}{i-1}}\right)^{-1}
((C\alpha)^{-[j]_q})^{-1}.
\end{eqnarray*}
En utilisant la formule du produit, le d\'eveloppement (\ref{ast}) et la formule \'ecrite plus haut avec $\mu=(1+qC\al)h$, on voit que la s\'erie $\dis\sum_{s\geq 0}\|\deform{q}{h}{L}_s\|x^{-\puiss{h}{q}{s}}$ admet pour s\'erie majorante la s\'erie qui constitue le d\'eveloppement en s\'erie de $(q,h)$-factorielles de la fonction : $$C\al\left(2+\frac1{x-\la}\right)\left(1+\frac{qC\al h}{x-\mu}\right).$$
Il s'ensuit que la famille 
$\deform{q}{h}{L}_{s}$ est fuchsienne $(C'',\lambda'')$ avec $\la''=\max(\la, (1+BC\al )B)$ si $B$ est un majorant de $\Lambda$.
Nous sommes donc en mesure d'appliquer la proposition \ref{preparatif} et cela termine la démonstration du théorème \ref{outil}.

\subsection{Fin de la preuve du th\'eor\`eme \ref{T1}}

La preuve du th\'eor\`eme \ref{T1} se fait maintenant de la manière suivante. On a d\'ej\`a remarqu\'e que la condition $c_i-c_j\not\in \rinteger$ assure la non r\'esonnance. Un calcul facile montre ensuite que la  matrice  $\deform{q}{h}{C_{0}}$ est conjugu\'ee \`a la matrice $diag\, (q^{c_1I_{r_1}+N_{r_1}},\cdots ,q^{c_mI_{r_m}+N_{r_m}})$ et  la premi\`ere hypoth\`ese du th\'eor\`eme \ref{outil} est satisfaite d\`es qu'il existe une famille born\'ee de conjugaisons. L'op\'erateur $\deform{q}{h}{\widetilde L_s}$ est l'inverse d'un op\'erateur qui  s'\'ecrit (toujours \`a conjugaison pr\`es) $U\mapsto diag\,(B_1,\cdots ,B_m)$ o\`u $B_j=q^{-c_j}[s]_qU+\tilde N_{r_j}U-q^sU\tilde N_{r_j}$ et la deuxi\`eme hypoth\`ese  d\'ecoule du fait   que $[s]_q$ et $q^s$ tendent vers l'infini avec $s$. Les autres affirmations du th\'eor\`eme sont  cons\'equences des r\'esultats du paragraphe pr\'ec\'edent et de la remarque  suivante. Si $\deform{q}{h}{F}(x)P{\cal E}_{\deform{q}{h}{J_{0}}}(x)$ est solution de (\ref{syst}), alors ${\deform{q}{h}{A}}^{\deform{q}{h}{F}(x)P}=\deform{q}{h}{J_{0}}$ et, puisque ${\deform{q}{h}{A}}^{\deform{q}{h}{F}(x)}=\deform{q}{h}{A_0}$, cela implique ${\deform{q}{h}{A_0}}^{P}=\deform{q}{h}{J_0}$, ce qui \'equivaut \`a $P\deform{q}{h}{J_{0}}=\deform{q}{h}{A_{0}}P$.

\section{Preuve du Th\'eor\`eme 2}

 Il s'agit ici d'\'etudier la continuit\'e en $(q,h)$ de la solution canonique obtenue au théorème \ref{T1} lorsqu'on suppose que la matrice $\deform{q}{h}{A}(x)$ est continue (en un sens \`a pr\'eciser) en un point $(q_0,h_0)$. Nous traitons ce probl\`eme lorsque $(q_0,h_0)$ est un point de la fronti\`ere de ${\cal Q}^+$, c'est-\`a-dire lorsque  $q_0=1$ ou $h_0=0$. Dans ce cas il est habituel de parler plut\^ot de confluence. 
 
  La forme de la solution canonique indique qu'on peut ici aussi s\'eparer le probl\`eme en deux parties que nous aborderons successivement : la confluence de la transformation de jauge d'une part et celle de la partie \og log - car \fg{} provenant de la matrice constante d'autre part.

\subsection{Confluence de la partie log - car}\label{convcl}
La famille des \og mon\^omes\fg{} $\left( x^{- \puiss{h}{q}{s}}\right)$, qui sont les caract\`eres correspondant \`a $c$ entier n\'egatif, est continue sur  ${\cal Q}^+$. C'est également le cas de la famille de polyn\^omes $(\deform{q}{h}{e_{c}}(x))$ o\`u $c\in\integer$.  Par contre, lorsque $c$ n'est pas entier, nous verrons que la famille des caract\`eres choisis et des logarithmes associ\'es n'a pas de limite en certains points du bord de ${\cal Q}^+$. On rencontrera cette situation dans le cas {\bf C1}; en revanche, on énoncera une propri\'et\'e de \og normalit\'e\fg{} et on déterminera l'ensemble des valeurs d'adhérence. L'absence de limite n'est pas un ph\'enom\`ene  surprenant car la d\'efinition de $\deform{q}{h}{e_{c}}(x)$ pour $q>1$, $h>0$ provient du cas $h=0$ par l'utilisation du changement de variable homographique  $t=(q-1)\frac xh+1$. Ce changement de variable n'a pas de limite (qui soit un changement de variable) en un point du bord de ${\cal Q}^+$ autre que $(1,0)$. Lorsque $(q_0,h_0)=(1,0)$, il faut supposer que $\frac {q-1}{h}$ a une limite si l'on veut que le changement de variable en ait une. 

Signalons un autre probl\`eme. L'\'equation $\sigma_qy(x)=q^cy(x)$ admet aussi la solution, in\-d\'e\-pen\-dan\-te de $q$, mais ramifi\'ee : $x^c$. La fonction $q$-p\'eriodique qui relie les deux caract\`eres s'exprime  \`a l'aide de la fonction $K(x,t)$ introduite par De Sole et Kac dans \cite{DSK}. Cette remarque est l'objet du lemme ci-dessous. Nous reprenons les notations suivantes de \cite{DSK}. Pour $a, t\in\complex$, on note :
 $$(1+a)_p^t=\frac{(-a;p)_\infty}{(-p^ta;p)_\infty}.$$ Si $t\in\complex$ et si on a choisi une d\'etermination   du logarithme de $x\neq 0$, on pose : $$K(x,t)=x^t\left(1+\frac1x\right)_p^t(1+px)_p^{-t}.$$
 Cette fonction v\'erifie les propri\'et\'es suivantes :
 \begin{enumerate}
 \item $K(x,t+1)=p^tK(x,t)$,
\item $K(qx,t)=K(x,t)$ si l'on suppose $\arg (qx)=\arg x$. Autrement dit, comme fonction de $x$, $K(x,t)$ est une $q$ - constante,
\item $K(x,t)=K(\frac1x,1-t)$ si l'on impose $\arg \frac1x=-\arg x$ de sorte que $\left(\frac1x\right)^t=x^{-t}$,
\item si $n\in\rinteger$, $K(x,n)=p^{\frac{n(n-1)}2}$ et est donc ind\'ependant de $x$.
\end{enumerate}

\begin{lem}\label{L100}
Pour $q>1$ :  $$\dis \deform{q}{0}{e_{c}}(x)=(-x)^cK(-\frac1x,c)=(-x)^cK(-x,1-c).$$
\end{lem}
 \begin{proof} Par d\'efinition :
  $$ \deform{q}{0}{e_{c}}(x) = \frac{\Theta_{q}(x)}{\Theta_{q,q^c}(x)}$$ et en utilisant la formule du triple produit  :
 \begin{eqnarray*}\deform{q}{0}{e_{c}}(x) &=& \frac{(x;p)_\infty( p/x;p)_\infty}{(p^cx;p)_\infty(p^{1-c}/x;p)_\infty}\\&=&(1-x)_p^c\left(1-\frac px\right)_p^{-c}\\&=&(-x)^c K(-\frac1x,c)=(-x)^cK(-x,1-c).\end{eqnarray*}
 \end{proof}

 Avec cette notation, pour $q>1$ et $h>0$ on a :
 $$\deform{q}{h}{e_{c}}(x) = \left(\frac{h}{1-p}\right)^{c}(1-z)_p^c$$ en posant $z=(1-p)\frac xh+p$.
 
  \subsubsection{Confluence des caractères}

Dans ce paragraphe, on s'int\'eresse au  probl\`eme suivant. Soit $(q_0,h_0)$ un point du bord de ${\cal Q}^+$. La fonction $\deform{q}{h}{e_{c}}(x)$ tend-elle vers la fonction $\deform{q_0}{h_0}{e_{c}}(x)$ lorsque $(q,h)$ tend vers $(q_0,h_0)$? Constatons tout d'abord que lorsque $q_0>1$ et $h_0=0$ la r\'eponse est n\'egative. Reprenons les notations pr\'ec\'edentes et d\'efinissons, pour $p=1/q$ fix\'e, le r\'eel $\xi$ par la condition
 $h=(1-p)p^\xi$ de sorte que $h\to 0$ si et seulement si $\xi\to+\infty$. On a :
 \begin{eqnarray*}
 \deform{q}{h}{e_{c}}(x)&=&p^{\xi c}(1-p^{-\xi}x-p)_p^c=p^{\xi c}\left(1-p^{-\xi}(x-p^{1+\xi})\right)_p^c\\
 &=&p^{\xi c}\frac{K(-\frac{p^\xi}{x-p^{1+\xi}}, c)}{\left(1-\frac{p^{\xi+1}}{x-p^{1+\xi}}\right)_p^{-c}}\,\left(-p^{-\xi}(x-p^{1+\xi})\right)^c\\&=&\frac{K(-\frac{p^\xi}{x-p^{1+\xi}}, c)}{\left(1-\frac{p^{\xi+1}}{x-p^{1+\xi}}\right)_p^{-c}}\,\left(-(x-p^{1+\xi})\right)^c\\&=&\frac{K(-\frac{p^{\{\xi\}}}{x-p^{1+\xi}}, c)}{\left(1-\frac{p^{\xi+1}}{x-p^{1+\xi}}\right)_p^{-c}}\,\left(-(x-p^{1+\xi})\right)^c
 \end{eqnarray*}
en d\'esignant par $\{\xi\}\in \mathop[0,1\mathop[ $ la partie fractionnaire de $\xi$. En utilisant le fait que $p^\xi\to 0$ quand $\xi\to +\infty$, on voit que le d\'enominateur tend vers $1$ et que le facteur $\left(-(x-p^{1+\xi})\right)^c$ tend vers le caract\`ere $(-x)^c$. Une condition n\'ecessaire \`a l'existence d'une limite est alors que $\xi\to+\infty$ de sorte que $\{\xi\}$ ait une limite. Or, m\^eme dans ce cas,  le facteur $K(-\frac{p^{-\{\xi\}}}x,c)$ est $q$-constant mais non constant. Il n'aurait donc pas \'et\'e plus int\'eressant  de faire un autre  choix de caract\`eres pour les op\'erateurs aux $q$-diff\'erences.

En utilisant  le lemme \ref{L100}, il vient  :

 $$\deform{q}{h}{e_{c}}(x)=\frac{\left(-p^{-\{\xi\}}(x-p^{1+\xi})\right)^{-c}}{\left(1-\frac{p^{\xi+1}}{x-p^{1+\xi}}\right)_p^{-c}}\,\left(-(x-p^{1+\xi})\right)^c \deform{q}{0}{e_{c}}(p^{-\{\xi\}}(x-p^{1+\xi}))$$ et donc :
 \begin{equation}\label{bbb}
 \deform{q}{h}{e_{c}}(x)=\frac{p^{\{\xi\}c}}{\left(1-\frac{p^{\xi+1}}{x-p^{1+\xi}}\right)_p^{-c}}\; \deform{q}{0}{e_{c}}(p^{-\{\xi\}}(x-p^{1+\xi})).
\end{equation}
 Cette expression montre qu'il n'y a pas de limite si $\xi\to+\infty$ mais que, si $\xi=\xi_0+n$ avec $\xi_0$ fix\'e dans $[0,1[$ (ou si, plus généralement, $(\xi_n)_{n \geq 1}$ est une suite telle que $\{\xi_n\}$ tend vers $\{\xi_0\}$ lorsque $n$ tend vers l'infini) il y a une limite quand $n\to+\infty$ qui vaut :
 $$p^{\xi_0c}\deform{q}{0}{e_{c}}(p^{-\xi_0}x)=p^{\xi_0c}\frac{\deform{q}{0}{e_{c-\xi_0}}(x)}{\deform{q}{0}{e_{-\xi_0}}(x)}.$$

Le m\^eme calcul montre aussi l'absence de limite dans le cas {\bf C3} sous l'hypoth\`ese {\bf H} avec $\ell=0$.

Ces propri\'et\'es \og n\'egatives\fg{} donnent de l'int\'er\^et au r\'esultat  de normalit\'e  qui figure dans le th\'eor\`eme suivant. Nous commen\c cons par \'enoncer un lemme  dont la preuve se trouve dans \cite{sauloy} (p. 1045).
\begin{lem} \label{lelem}Soient $p$ un r\'eel tel que $0<p<1$, $a \in \complex$, $b\in \complex \setminus p^{-\integer}$. On pose  $d=|a-b|$ et $\dis \eta = \inf_{s \in \integer} |1-p^{s}b|$. Alors, pour tout $s \in \integer$ :
  \begin{equation*}
  \left|\frac{(a;p)_{s}}{(b;p)_{s}}\right|
  \leq\exp(\frac{d}{\eta}\frac{1}{1-p}).
  \end{equation*}
\end{lem}

\begin{theo}\label{Convc} $_{}$ \vskip 1 pt
\begin{description}
\item[C1]  Pour tout $\ve>0$, la famille $\left(\deform{q_0}{h}{e_{c}}(x)\right)_{0<h\leq\ve}$ est normale sur $\complex \setminus q_0^{c}\real^+$. 

Si  $\xi\in\mathop[0,1\mathop[$ la suite   $\left(\deform{q_0}{h_n}{e_{c}}\right)_{n\geq 0}$ o\`u $h_n=(1-p_0)p_0^{\xi +n}$ ($p_0=1/q_0$) tend vers
$$\tilde e_c^{(q_0,\xi)}(x)=p_0^{\xi c}\frac{\Theta_{q_0,p_0^{\xi}}}{\Theta_{q_0,p_0^{\xi}p_0^{-c}}}=p_0^{\xi c}\frac{\deform{q_0}{0}{e_{c-\xi}}(x)}{\deform{q_0}{0}{e_{-\xi}}(x)}$$ uniform\'ement sur tout compact de $\complex \setminus q_0^{c}\real^+ $. R\'eciproquement si la suite $\left(\deform{q_0}{h_n}{e_{c}}(x)\right)_{n\geq 0}$ converge lorsque $h_n\to 0$, il existe $\xi\in\mathop[0,1\mathop[$ tel que sa limite soit  $\tilde e_c^{(q_0,\xi)}(x)$.
\item[C2] Pour tout $c\in\complex$, la famille $\left(\deform{q}{h}{e_{c}}\right)_{(q,h)\in{\cal Q}}$ converge 
 uniform\'emement sur tout compact de $\complex\setminus(c+h_0\integer^*)$ vers $\deform{1}{h_0}{e_{c}} $.
\item[C3  ] On suppose $\ell\neq 0$.  Pour tout $c\in\complex$, la famille $\left(\deform{q}{h}{e_{c}}\right)_{(q,h)\in{\cal Q}^*}$ converge 
 uniform\'emement sur tout compact de $\complex \backslash \real^{+}$ vers $e_{c}$.
\end{description}
\end{theo}
\begin{proof}$_{}$\vskip 1 pt
\begin{description}
\item[C1]
Rappelons l'égalité :

\begin{eqnarray*}
\deform{q_0}{h}{e_{c}}(x) &=& \frac{K(-\frac{p_0^{\{\xi\}}}{x-p_0^{1+\xi}}, c)}{\left(1-\frac{p_0^{\xi+1}}{x-p_0^{1+\xi}}\right)_{p_0}^{-c}}\,\left(-(x-p_0^{1+\xi})\right)^c
\end{eqnarray*}
où le r\'eel $\xi$ est défini par
 $h=(1-p_0)p_0^\xi$.
Soit $K$ un compact de $\complex \setminus q_0^{c}\real^+$.
Puisque $\xi\to+\infty$ lorsque $h\to 0$, il existe $A_1>0$
tel que, pour tout $\xi \geq A_1$ et pour tout $x \in K$,
$\left|\frac{p_0 ^{\xi+1}}{x-p_0^{1+\xi}}\right|\leq 1/2 .$ On en d\'eduit qu'il existe $M_1>0$ tel que pour tout $\xi \geq A_1$
et pour tout $x \in K$,
$\frac{1}{\left|\left(1-\frac{p_0 ^{\xi+1}}{x-p_0^{1+\xi}}\right)_{p_0}^{-c}\right|} \leq M_1.$
D'autre part, puisque $K$ est un compact de $\complex \setminus q_0^{c}\real^{+}$ et que
$p_0^\xi$ tend vers $0$ avec $h$, il existe $A_2>0$ et $\theta\in\mathop]0,\pi\mathop[$
tels que pour tout $\xi \geq A_2$ et pour tout $x \in  K$,
$\left|arg(\frac{p_0^{\{\xi\}}}{x-p_0^{1+\xi}}q_0^c)\right|\geq \theta$.
Ainsi, puisque $p_0 \leq p_0^{\{\xi\}} \leq 1$, il existe un compact $K'$
de $\complex \setminus q_0^{c}\real^+$ contenant
$\frac{p_0^{\{\xi\}}}{x-p_0^{1+\xi}}$ pour tout $\xi \geq A_2$
et pour tout $x \in  K$. Il existe donc $M_2>0$ tel que pour tout $\xi \geq A_2$
et tout $x \in K$, $\left|K(-\frac{p_0^{\{\xi\}}}{x-p_0^{1+\xi}}, c)\right|\leq M_2.$
On en déduit que la famille des fonctions $\deform{q_0}{h}{e_{c}}$, $h>0$ petit est
uniformément bornée sur $K$. Cela prouve la normalité.

Les deux autres assertions sont des conséquences immédiates de la formule (\ref{bbb}). Pour la réciproque, il suffit d'utiliser le fait que la
suite des
$\{\xi_n\}$ admet toujours au moins une valeur d'adhérence.

\item[C2] Le cas  $h=h_0=1$ est trait\'e dans \cite{duv1}. On emploie la m\^eme m\'ethode. On fait le changement de variable 
 : $x=h[u]_q$ qui permet d'\'ecrire $\deform{q}{h}{e_{c}} $ \`a l'aide de la fonction $\Gamma_p$ de Jackson d\'efinie par :
$$\Gamma_p(t)=\frac{(1-t)_p^{t-1}}{(1-p)^{t-1}}=(1-p)^{1-t}\frac{(p;p)_\infty}{(p^t;p)_\infty}.$$
En effet, en remarquant que : $$z=(1-p)\frac xh+p=(1-p)\frac{q^u-1}{q-1}+p=p(p^{-u}-1)+p=p^{1-u},$$ on a :
$$\deform{q}{h}{e_{c}} (x)=h^c\frac{\Gamma_p(1+c-u)}{\Gamma_p(1-u)}.$$ En utilisant le d\'eveloppement en produit infini de la fonction  $\Gamma_p$ dans le demi-plan $\Re z>0$ :
$$\Gamma_p(z)=\prod_{n\geq 1}\frac{(1-p^{n+1})^z}{(1-p^{n+z})(1-p^n)^{z-1}}$$ on obtient le d\'eveloppement :
$$\deform{q}{h}{e_{c}} (x)=h^c\prod_{n\geq 1}\left(\frac{1-p^{n+1}}{1-p^n}\right)^c\frac{1-p^{n-u}}{1-p^{n+c-u}}.$$ Le terme g\'en\'eral de ce produit s'\'ecrit aussi :
$$\left(\frac{q^{n+1}-1}{q^n-1}\right)^c\frac{q^n-1-(q-1)\frac xh}{q^{n+c}-1-(q-1)\frac xh}=\left(\frac{[n+1]_q}{[n]_q}\right)^c\frac{[n]_q-\frac xh}{[n+c]_q-\frac xh}$$ et tend, quand $q\to 1^+$, vers $\dis\left(\frac{n+1}n\right)^c\frac{n-\frac xh}{n+c-\frac xh}$, terme g\'en\'eral de l'expression en produit infini de $\dis\frac{\Gamma(1+c-\frac xh)}{\Gamma(1-\frac xh)}$. Les arguments de convergence sont alors ceux de \cite{koo} par exemple.

\item[C3] La convergence unifome sur tout compact du disque unité, lorsque $q$ tend vers $1$ par valeurs r\'eelles positives, de $z \mapsto \frac{(z;p)_{\infty}}{(p^{c}z;p)_{\infty}}$ vers $z\mapsto (1-z)^c$ est une conséquence classique du théorème $q$-binomial (où on a utilisé la détermination principale du logarithme). De plus, la famille des $z \mapsto \frac{(z;p)_{\infty}}{(p^{c}z;p)_{\infty}}$, $p<1$ est normale sur $\complex \setminus [1,+\infty[$ : cela résulte du lemme \ref{lelem}. On en déduit que $z \mapsto \frac{(z;p)_{\infty}}{(p^{c}z;p)_{\infty}}$ converge uniformément sur tout compact de $\complex \setminus [1,+\infty[$ vers $z\mapsto (1-z)^c$ lorsque $q\rightarrow 1^+$ et le résultat suit dans le cas o\`u $\dis\frac h{q-1}\to \ell\neq 0$.

\end{description}
\end{proof}
On peut  calculer  la moyenne des
caractères obtenus par passage à la limite dans le cas \textbf{C1}
(en les considérant comme paramétrés par $\alpha=p_0^\xi$). Nous
pouvons envisager deux types de moyennes : celle obtenue grâce à
une $q$-intégrale et celle calculée à l'aide d'une intégrale
classique (mesure de Lebesgue). Rappelons que, si $a$ et $b$ sont
deux nombres réels positifs avec $b>a$ et si $f$ est une fonction
à valeurs complexes définie sur $[0,b]$, alors la $q$-intégrale de
$f$ entre $a$ et $b$ est donnée par :
\begin{eqnarray*}
\int_a^b f (\alpha) d_p \alpha &=&  \int_0^b f (\alpha) d_p \alpha - \int_0^a f (\alpha) d_p \alpha \\
&=& (1-p)\sum_{s=0}^{+\infty} f(p^s b) p^s - (1-p)\sum_{s=0}^{+\infty} f(p^s a) p^s.
\end{eqnarray*}
Dans notre cas un calcul élémentaire donne la valeur explicite suivante pour la $q$-moyenne :  $$\int_{p_0}^1 \alpha^c \frac{\Theta_{q_0,\alpha}}{\Theta_{q_0,\alpha
    p_0^{-c}}}(x)d_{p_0}\alpha= (1-p_0) \frac{\Theta_{q_0}}{\Theta_{q_0,p_0^{-c}}}(x)=(1-p_0)e_c^{(q_0,0)}(x).$$
En ce qui concerne la seconde moyenne, qu'on note $\phi$ :
$$
\phi(x)=\int_p^1 \alpha^c \frac{\Theta_{q,\alpha}}{\Theta_{q,\alpha
    p^{-c}}}(x) d \alpha=(-x)^c\int_p^1K(-tx^{-1},c)\,dt$$
nous ne savons pas en faire un calcul explicite. Soulignons cependant qu'elle satisfait l'équation différentielle suivante qui la relie au r\'esultat pr\'ec\'edent:

\begin{eqnarray*}
\phi'(x)&=& \frac{c+1}{x} \phi(x) - \frac{1-p}{x}e_c^{(q,0)}(x)
\end{eqnarray*}
et qu'elle n'est pas définie sur $\complex^*$ ni même sur son revêtement universel mais seulement sur le revêtement universel de $\complex^* \setminus q^{c+\rinteger}.$
\subsubsection{Confluence des logarithmes}

En remarquant que $\deform{q}{h}{e_{c}}(x)$ est holomorphe en $c$, on peut \'echanger les d\'erivations successives en $c$ et les op\'erations de passage \`a la limite. Reprenant alors les divers arguments de la preuve du th\'eor\`eme \ref{Convc} on \'etablit la proposition suivante.

\begin{prop}\label{Convl} Soit $k$ un entier strictement positif.
\begin{description}
\item[C1]  Pour tout $\ve>0$, la famille $\left(\deform{q_0}{h}{\ell_{c,k}}(x)\right)_{0<h\leq\ve}$ est normale sur $\complex \setminus q_0^{c}\real^+$. La suite $\left(\deform{q_0}{h_n}{\ell_{c,k}}\right)_{n\geq 0}$ o\`u $h_n=(1-p_0)p_0^{n}$  tend vers $\deform{q_0}{0}{\ell_{c,k}}$, uniform\'ement sur tout compact de $\complex \setminus q_0^{c}\real^+$.
\item[C2] La famille $\left(\deform{q}{h}{\ell_{c,k}} \right)_{(q,h)\in{\cal Q}}$ converge uniform\'emement sur tout compact de $\complex\setminus(c+h_0\integer^*)$ vers $\deform{1}{h_0}{\ell_{c,k}} $.
\item[C3  ] Si $\ell\neq 0$, la famille $\left(\deform{q}{h}{\ell_{c,k}}\right)_{(q,h)\in{\cal Q}}$ converge uniform\'ement sur tout compact de $\complex \backslash \real^{+}$ vers $\deform{1}{0}{\ell_{c,k}}$ .
\end{description}
\end{prop}

\subsection{Confluence de la transformation de jauge}

Pour \'etablir le th\'eor\`eme \ref{T2}, il reste à étudier la partie transformation de jauge. On remarque d'abord qu'en  passant \`a la limite dans les relations $\|\deform{q}{h}{A}_s\|\leq Cq^{s-1}\left(\la^{-\puiss{h}{q}{s-1}}\right)^{-1}$ données par les hypothèses du th\'eor\`eme \ref{T2} et en utilisant la continuit\'e en $(q,h)\in{\cal Q}^+$ du majorant, on obtient la premi\`ere affirmation du théorème gr\^ace \`a la proposition \ref{prop1}. L'hypoth\`ese que le syst\`eme aux $(q_0,h_0)$-diff\'erences de matrice ${A}(x)$ est fuchsien signifie que la matrice $C_0=I-(q_0-1)A_0$ appartient \`a $Gl_n(\complex)$. Si ce syst\`eme est non r\'esonnant, on peut appliquer le th\'eor\`eme \ref{T1}  \`a la famille $\Lambda'=\Lambda\cup\{(q_0,h_0)\}$ en posant $\deform{q_0}{h_0}{A}(x)=A(x)$. Il existe donc une transformation de jauge $\deform{q_0}{h_0}{F}(x)$ telle que $A^{\deform{q_0}{h_0}{F}}=A_0$. La preuve du th\'eor\`eme \ref{T2} est alors cons\'equence du r\'esultat suivant qui \'etablit la continuit\'e de la transformation de jauge.  

\begin{theo} Sous les hypoth\`eses du th\'eor\`eme \ref{T2},  l'unique transformation de jauge $\deform{q}{h}{F}(x)$, $(q,h) \in \Lambda$ donn\'ee par le th\'eor\`eme \ref{T1} converge  vers l'unique transformation de jauge $\deform{q_0}{h_0}{F}(x)$ tangente \`a l'identit\'e qui conjugue le syst\`eme aux $(q_0,h_0)$-diff\'erences de matrice ${A}(x)$ \`a la matrice constante $A_0$.
La convergence est uniforme sur tout compact de $V^{(q_0,h_0)}(\la')$ dans les cas {\em {\bf C1}} et {\em {\bf C2}}, de $V_\ell(\la')$ dans le cas {\em {\bf C3}}.
 \end{theo}

 \begin{proof} Le th\'eor\`eme \ref{T1} donne des constantes $C',\la'>0$ telles que, si $\dis \deform{q}{h}{F}(x)=I+\sum_{s\geq 1}\deform{q}{h}{F}_sx^{-\puiss{h}{q}{s}}$, on ait, pour tout $(q,h)\in\Lambda$ :
 $$ \|\deform{q}{h}{F}_s\|\leq C'q^{s-1}\left(\la'^{-\puiss{h}{q}{s-1}}\right)^{-1}\leq \frac{C'}{\la'}\left(\la'^{-\puiss{h}{q}{s}}\right)^{-1}$$ et donc:
 $$\|\deform{q}{h}{F}_sx^{-\puiss{h}{q}{s}}\|\leq \frac{C'}{\la'}\prod_{k=0}^{s-1}\frac{q^k\la'+h[k]_q}{|q^kx+h[k]_q|}.$$ Cette majoration permet de  traiter chaque cas comme suit :
 \begin{enumerate}
  \item[cas {\bf C2}] Le
  changement de variable $x \leftarrow \frac h{h_0}x$ permet   d'appliquer le th\'eor\`eme  3.1 de \cite{duv1}.
  \item[cas {\bf C1}]
  Soient $q_1\geq q_0$ un majorant de $\{q\;|\;(q,\la)\in \Lambda\}$ et $R>\la'$. Pour tout $x$ tel que $|x|\geq R$ et tout entier $k\geq 0$,  l'inégalité triangulaire renversée jointe à l'inégalité $|-\puiss{h}{q}{k}q^{-k}| \leq \frac{h}{q-1}$, implique que $|x+hq^{-k}[k]_q|\geq R-\frac h{q_1-1}$. Supposons que $\ve>0$ est assez petit pour que $R-\frac h{q_1-1}>0$ si $h<\ve$. Alors, en posant $\kappa=\frac1{q_1-1}$, on a :
  $$\frac{q^k\la'+h[k]_q}{|q^kx+h[k]_q|}=\frac{\la'+hq^{-k}[k]_q}{|x+hq^{-k}[k]_q|}\leq \frac{\la'+h\kappa}{R-h\kappa}\leq \frac{\la'+\ve\kappa}{R-\ve\kappa}=\rho<1$$ si  $\ve$ est assez petit. Ainsi, on peut majorer le terme g\'en\'eral de la s\'erie $\deform{q}{h}{F}(x)$ par $M\rho^s$ et le r\'esultat s'en d\'eduit par convergence dominée.

 \item[cas {\bf C3}] Le cas $\ell\in \real^+$ se traite comme le cas {\bf C1} en rempla\c cant $\la'$ par $\la'-\ell$ (suppos\'e positif). Le cas $\ell=+\infty$ se traite quant à lui comme le cas {\bf C2}. En effet, on a l'estimation suivante : $$
   \left|\frac{q^{k}\lambda+\puiss{h}{q}{k}}{q^{k}x+\puiss{h}{q}{k}}\right|
   \leq
   \frac{q^{k}\lambda+\puiss{}{q}{k}}{q^{k}\Re x+\puiss{}{q}{k}}.
  $$
et on remarque que le majorant obtenu est indépendant  de $h$ : on termine maintenant la preuve par les mêmes raisonnements que pour le cas {\bf C2}. 
\end{enumerate}
 \end{proof}

Signalons les points suivants.  Le cas trait\'e dans \cite{duv1} est le cas {\bf C2} avec $h=1$. Le cas trait\'e dans \cite{roques} est le cas {\bf C3} avec $q=1$ et donc $\ell=+\infty$. Le cas trait\'e dans \cite{sauloy} est le cas {\bf C3} avec  $h=0$ et donc $\ell=0$. 

\bigskip

Notons pour terminer que le th\'eor\`eme \ref{T2} s'applique en particulier  lorsqu'on d\'eforme le syst\`eme non r\'esonnant de matrice :

$$A(x)=\sum_{s\geq 0}A_sx^{-\puiss{h_0}{q_0}{s}}$$ o\`u la s\'erie converge dans $V^{(q_0,h_0)}(\la)$ en la famille de syst\`emes de matrice:
$$A^{(q,h)}(x)=A_0^{(q,h)}+\sum_{s\geq 1}A_sx^{-\puiss{h}{q}{s}}$$ o\`u, si: $$A_0=P \,diag \,(\deform{q_0}{h_0}{A_{c_1,r_1}},\cdots, \deform{q_0}{h_0}{A_{c_m,r_m}})\,P^{-1},$$ alors:

$$A_0^{(q,h)}=P \,diag \,(\deform{q}{h}{A_{c_1,r_1}},\cdots, \deform{q}{h}{A_{c_m,r_m}})\,P^{-1}.$$

\end{document}